\documentclass[11pt]{article}
\usepackage{graphicx,epsfig,amssymb,cite}
\usepackage{amsmath}

\usepackage{algorithm}
\usepackage{algpseudocode}
\usepackage{mathrsfs}
\usepackage[utf8]{inputenc}
\usepackage{amsmath,amssymb}  

\newtheorem{theorem}{Theorem}
\newtheorem{remark}{Remark}
\newtheorem{lemma}{Lemma}

\newtheorem{assumption}{Assumption}

\newtheorem{proposition}{Proposition}
\newtheorem{definition}{Definition}

\date{}

\begin{document}
\title{Learning-based primal-dual  optimal control of discrete-time stochastic systems with multiplicative noise
}

\author{
Xiushan Jiang$^1$, Weihai Zhang$^2$\thanks{Corresponding author: Weihai Zhang (email: w\_hzhang@163.com).} \\
$^1$ \small College of New Energy, China University of Petroleum (East China), Qingdao 266580, China \\
$^2$ \small College of Electrical Engineering, Shandong University of Science and
\\ \small Technology, Qingdao 266590, China
}
\maketitle

{\bf Abstract-}Reinforcement learning (RL) is an effective approach for solving
optimal control problems without knowing the exact information of
the system model. However, the classical  Q-learning method, a
model-free RL algorithm, has its limitations, such as lack of strict
theoretical analysis and the need for artificial disturbances during
implementation. This paper explores the partially  model-free stochastic linear quadratic  regular
 (SLQR)  problem  for a system with multiplicative noise from the primal-dual perspective to address these challenges.
This
approach lays a strong theoretical foundation for understanding the
intrinsic mechanisms of classical  RL algorithms.  We reformulate
the SLQR into  a  non-convex  primal-dual  optimization  problem and  derive a
strong duality result, which enables  us to provide  model-based and model-free algorithms for SLQR  optimal policy design
based on the Karush-Kuhn-Tucker (KKT) conditions. An illustrative example
demonstrates the proposed model-free algorithm's validity,
showcasing the central nervous system's learning mechanism in human
arm movement.


\section{Introduction}

Linear quadratic regular (LQR) was initiated by Kalman
\cite{kalman}, and further developed in \cite{anderson,lewis}. It is
well-known  that    LQR  is one of the most important
optimal controls, which is very elegant in theory and has more
applications in engineering practice \cite{dom,zhanghan1}.
Stochastic linear quadratic regular (SLQR) seems to be first studied
by \cite{wonham}, in particular, since \cite{che_98} established the
indefinite SLQR theory, SLQR has gained a lot of scholars'
attention, and has been extensively studied; see
\cite{huangyulin,sun,ydd,zhang-xie-chen}.
Generally speaking, SLQR will lead to solving a  generalized
algebraic Riccati equation (GARE), which requires us to know
complete system information including the system structure and exact
parameter information. However, the exact model structure and system parameters are commonly unknown in the process of  practical modeling, in this case, all traditional model-based methods  become invalid.
The model-free reinforcement learning (RL) approach provides a solution to unknown dynamics  by exploring poorly structured systems through
state-input data analysis.

RL is a branch  of machine learning that iteratively achieves an optimal policy through interactions with
the environment. Pioneering studies  \cite{rlcon2,rlcon1}  initiated
the development of RL within the optimal control framework, which has since garnered significant attention \cite{Busoniu,pangbo1,oura2024,youbased,wang2025,jiangsci}. For example, the reference  \cite{jiangzp} used the off-policy RL to study the $H_\infty$ control of linear discrete-time systems, while
\cite{vrabie} researched the adaptive optimal control of linear continuous-time systems based on the  policy iteration (PI). In  \cite{pangbo}, the authors used a novel off-policy RL method named optimistic
least squares-based PI  to find
directly  near-optimal controllers  from input/state data
 for adaptive optimal stationary  control of
linear It\^o  systems with additive and multiplicative noises.
Q-learning is one of the most important RL approaches, which has
 been  studied  as an effective model-free RL algorithm \cite{farja,Fazel,Karl,shuzhan,Stephen,ybbasi} for  LQR optimal policy design. This is because LQR is  one class of the most important and simplest  optimal controls, which captures the main characteristics of  Q-learning.  Particularly, in \cite{Fazel}, a stochastic policy gradient algorithm was presented  with a shortcoming of large variance. Modified  approximate  PIs for LQR were given in \cite{Karl,ybbasi}.  In \cite{shuzhan}, the model-free RL algorithm based on Q-function was  proposed for discrete-time systems with multiplicative and additive  noises.  The authors of \cite{Stephen} discussed the gap between model-based  and model-free model algorithms  on LQR. In our recent work \cite{jiangauto}, we  applied an off-policy RL method to study the stochastic $H_\infty$ control problem with unknown system model.

As said in \cite{watkins},  Q-function learning-based  optimal policy
is only guaranteed for finite Markov decision process, and this  limitation poses significant data storage
requirements for complex systems  and makes  its practical
applications face  challenges.  It can be found that,  most of Q-learning
algorithms  lack of a solid theoretical analysis or scalability  and are  dependent  on
the persistent excitation  assumption.  To address these shortcomings, a novel primal-dual
  Q-learning framework for the LQR problem was recently
  introduced  \cite{lee2019} and further developed \cite{liman} to investigate the SLQR  problem
  with additive Gaussian white  noises.
In \cite{liman}, the random variables  are assumed to be Gaussian white noises, and the cost functions are quadratic with a discount factor $\gamma$ belonging to $0<\gamma<1$.
  Primal-dual
  RL   method clarifies  the  essential   relations  among Q-learning algorithms, off-line PI algorithm and LQR optimization based on
  semidefinite programming \cite{ydd}.  Primal-dual  RL  algorithms \cite{lee2019,liman,lixiuxian} possess the advantages of a fast convergence  and convenience for handling higher  dimensional systems.
  However,
  to the best of our knowledge, up to now, no literature has succeeded in solving the model-free optimal  policy design for the  SLQR of discrete-time multiplicative noise systems from the primal-dual perspective.

This paper aims to  explores the partially model-free RL
algorithms of SLQR in linear  discrete-time stochastic systems with
multiplicative noises  by  employing the primal-dual approach.  This paper can be viewed as a  non-trivial  extension of \cite{lee2019} to stochastic multiplicative noise systems  due to that   there have essential differences between deterministic and stochastic systems. In fact, in order to develop a parallel frame to deterministic primal-dual model-free algorithms, we have to apply our previously introduced new  definitions and theorems  such as exact observability, exact detectability, Popov-Belevith
-Hautus  (PBH) criteria for eigenvector
test of exact observability and exact detectability,  as well as generalized Lyapunov theorems for the asymptotical   stability  in mean-square (ASMS) sense \cite{huangyulin,zhang-xie-chen}, while related definitions and results of continuous-time It\^o systems can be found in \cite{zwh2004,zhang2008,zhangphd}.

The main contributions of this paper are as follows:
 \begin{description}

   \item[(1)]
  We propose a novel off-line PI to solve the GARE from the concerned SLQR, and  a strict convergence proof is also presented.
  More importantly, we point out that this PI algorithm has quadratic convergence speed; see Remark~\ref{eqvhvhaa}.
  Our off-line PI algorithm  can be viewed as a discretized  version of \cite{zhangphd}, which can also be viewed as an
  extension of classical Kleinman iteration algorithm \cite{klein}.

   \item[(2)]
 When the drift term coefficients  are  unknown,  we propose a novel primal-dual optimization algorithm  to obtain a partially
 model-free SLQR optimal policy
  design.  As corollaries, all  results of the deterministic system
   \cite{lee2019} can be obtained. When the diffusion term
   coefficients are also unknown, the fully primal-dual model-free
   SLQR optimal policy design remails unsolved.

   \item[(3)]
  Compared with the  persistent excitation condition-based Q-learning algorithm \cite{pangbo}, our primal-dual-based algorithm can
   quickly converge
  to the optimal solution. Moreover, the designed
  algorithm is obtained by solving the  Karush-Kuhn-Tucker (KKT) condition, which not only demonstrates the  equivalence with respect to  classical
   PI and
  Q-learning algorithms \cite{shuzhan}, but also provides a rigorous convergence analysis  for RL design of SLQR.
  \end{description}

The organization of  this paper is as follows:   In Section 2.1, we
first formulate  the SLQR problem, and make some preliminaries such
as exact observability/exact detectability, generalized Lyapunov
theorem  and PBH  criteria.  Then,  in Section 2.2,   we propose
an off-line   PI algorithm to solve GARE  with a  strict convergence
analysis and Q-learning function. In Section 3, we reformulate the
SLQR optimality  into a nonlinear constrained optimization problem
via constructing  a proper Lagrangian dual function, and then  prove
the strong duality which  yields the KKT condition. Based on KKT
condition, both model-based and partially model-free primal-dual
algorithms for searching for optimal control policy are given. An
illustrative example in Section 4 demonstrates the efficiency of the
proposed partially model-free algorithm.  Section 5 concludes this
paper with some remarks and future perspective.

  Notations: $\mathcal{C}$: the complex plane; $\mathscr{S}_n$: the collection of all $n\times n$  symmetric matrices;
  $\mathscr{S}_n^+$ ($\mathscr{S}_n^{++}$): the set of all $n\times n$ real symmetric positive semidefinite
  (positive definite) matrices;
  $\mathcal{N}_+$($\mathcal{N}$):  set of positive (non-negative) integers; ${\mathcal N}_T:=\{0,1,\dots, T\}$;
  $\|\cdot\|$:  the Euclidean vector norm or Frobenius matrix
norm;
  $P\succ 0$($\succeq 0$): $P$ is a positive definite
   (positive semidefinite) symmetric
    matrix; $\sigma(\mathcal{L})$: the spectrum set of the operator $\mathcal{L}$;  $ {\mathcal
D}(0,1):=\{\lambda\in {\mathcal C}: |\lambda|<1\}$;  $A'$: the
transpose of the matrix $A$; ${\mathcal L}^2_{{\mathcal F}_{k}}(\Omega, X)$: the family of $X$-valued
${\mathcal F}_{k}$-measurable  random   variables with bounded variances, i.e., for any
 $\xi$ from the family, ${\mathcal E}\|\xi\|^2<\infty$; $l^2_w({\mathcal N},{\mathcal R}^k)$:  the set of all
non-anticipative square summable stochastic processes
$$
u=\{u_k: u_k\in {\mathcal L}^2_{{\mathcal F}_{k-1}}(\Omega,
{\mathcal R}^{m})\}_{k\in {\mathcal N}}
$$
with the $l^2$-norm of $u\in l^2_w$
defined by
$$
\parallel u \parallel_{l^2_w}=\left(\sum^{\infty}_{k=0}{\mathcal E}
\parallel u_k \parallel^2\right)^{\frac{1}{2}}.
 $$

\section{Problem formulation and preliminaries}
\subsection{SLQR problem}

In this subsection, we consider the following  linear discrete-time
stochastic system with multiplicative noise
\begin{eqnarray}\label{sys1}
x_{k+1}=Ax_k+Bu_k+(Cx_k+Du_k)w_k,\ x_0=z\in\mathcal{R}^n,
\end{eqnarray}
where $x_k\in \mathcal{R}^n$ and $u_k\in\mathcal{R}^m$ with
$k\in\mathcal{N}$ are the state vector and control action,
respectively. $A, B, C$,  and $D$ are system matrices with suitable
dimensions. $\{w_k, k\in\mathcal{N}\}$ is a sequence of real
independent  random variables with $\mathcal{E}(w_k)=0$ and
$\mathcal{E}(w_kw_s)=\delta_{ks}$ (Kronecker function) which is
 defined over  a complete filtered
 probability space $\{\Omega, \mathcal{F},
 \mathcal{P}; \mathcal{F}_k\}$ with
  $\mathcal{F}_k$ being the $\sigma$-algebra generated by  {$\{w_v, v=0, 1, 2, ..., k-1\}$}.
Without loss of generality, we assume that the initial condition $x_0=z$ is
a deterministic vector.
For simplicity, the notation $[A, B; C, D]$  refers to the
system (\ref{sys1}).
\begin{remark}
All results of this paper can be generalized to multiple
multiplicative noise cases, i.e., system (\ref{sys1}) can be
generalized to
\begin{eqnarray*}
x_{k+1}=Ax_k+Bu_k+\sum_{i=1}^N(C_ix_k+D_iu_k)w^i_k,\
x_0\in\mathcal{R}^n,
\end{eqnarray*}
where $\{w_k^1\}_{k\in{\mathcal N}}$, $\cdots$,
$\{w_k^N\}_{k\in{\mathcal N}}$  are mutually independent random
variable
 sequences. Here, we consider system (\ref{sys1}) only for
simplicity.  In addition, we do not require $w_k$, $k\in {\mathcal N}$, obey the Gaussian distribution as done in \cite{shuzhan,liman}.
\end{remark}

The cost function associated with the system (\ref{sys1}) is denoted by
\begin{eqnarray}\label{cost1}
J(z, u)=\mathcal{E}\left[\sum_{k=0}^\infty(x_k'Qx_k+u_k'Ru_k)\right],
\end{eqnarray}
where Q and $R$ are symmetric matrices with appropriate dimensions with $Q\succeq 0$ and $R\succ0$.
\begin{definition}
System (\ref{sys1}) or  $[A,B; C,D]$ is called stabilizable  if
there exists  a  feedback  control policy $u_k=Fx_k$ with the
constant matrix $F$, such that for any  initial state $x_0=z$,   the
closed-loop system
\begin{eqnarray}\label{sysF}
x_{k+1}=(A+BF)x_k+(C+DF)x_kw_k
\end{eqnarray}
is ASMS, that is, we have
   $\lim_{k\rightarrow{+\infty}}\mathbb{E}[(x_k^{F, z})'(x_k^{F, z})]=0$, where
   the solution $x(k; F, z)$ of system (\ref{sysF}) is denoted as $x_k^{F, z}$ for simplicity.
  When (\ref{sysF}) is ASMS,  we also call  $[A+BF; C+DF]$  ASMS  for short.  Moreover, the feedback gain $F$ is called a stabilizing state-feedback gain.
\end{definition}
Under the state-feedback  gain $F$, the cost function (\ref{cost1}) is denoted by
\begin{eqnarray}\label{cost2}
J(z, F)=\mathcal{E}\left\{\sum_{k=0}^\infty\left[\begin{array}{ccc}x_k^{F, z}\\Fx_k^{F, z}\end{array}\right]'
    \Lambda\left[\begin{array}{ccc}x_k^{F, z}\\Fx_k^{F, z}\end{array}\right] \right\}
\end{eqnarray}
with $\Lambda =\left[\begin{array}{ccc}Q&0\\0&R\end{array}\right]$.

The SLQR problem can be stated as follows: Under the system (\ref{sys1}), search for, if it exists, an admissible control $u^*_k=F^*x_k\in {\mathcal U}_{ad}$  to minimize $J(z, F)$, where
\begin{eqnarray*}
&&{\mathcal U}_{ad}:=\{u\in l^2_w({\mathcal N}, {\mathcal R}^m): \{u_k\}_{k\in
{\mathcal N}} \ \ \mbox{is\ a \ mean\ square }\\
&&  \ \ \ \ \ \ \ \ \ \ \ \ \  \mbox{stabilizing\ control\ sequence}\}.
\end{eqnarray*}
In this case,  $\{u^*_k\}_{k\in {\mathcal N}}$ is  called the optimal control sequence, while
$\{x^*_k\}_{k\in {\mathcal N}}$  corresponding to $\{u^*_{k}\}_{k\in
{\mathcal N}}$ is the optimal state  trajectory, and $J(z, F^*)$ is the optimal
cost value.

\begin{definition}\label{def:Definition 3.2.3}
The system
\begin{equation}
\left\{
\begin{array}{l}
x_{k+1}=Ax_k+Cx_kw_k,\ \ x_0\in {\mathcal R}^n,\\
 y_k=Qx_k,  ~ k\in {\mathcal N}
\end{array}
\right. \label {eq ghvhvyg3.3.1}
\end{equation}
or $(A, C|Q)$ is said to be exactly observable  if there exists
$T\in {\mathcal N}_+$ such that
$$
y_k\equiv 0, a.s., \forall k\in {\mathcal N}_T \Rightarrow x_0=0.
$$
$(A, C|Q)$ is said to be exactly detectable if
$$
y_k\equiv 0, a.s., \forall k\in {\mathcal N}_T \Rightarrow
\lim_{k\to \infty}{\mathcal E}\|x_k\|^2=0.
$$
\end{definition}
\begin{remark}
Definition~\ref{def:Definition 3.2.3} loosens the conditions of
Definition 3.7 of \cite{zhang-xie-chen}, where $\forall k\in
{\mathcal N}$ in \cite{zhang-xie-chen} is replaced by $\forall k\in
{\mathcal N}_T$ for $T\in {\mathcal N}_+$.
\end{remark}

The following lemma can be found in Theorem 3.6  and  Lemma 3.5 of
\cite{zhang-xie-chen}.
\begin{lemma}\label{lem2.1}
    For the  system $[A,0; C,0]$ or $[A;C]$,
    the following three statements are equivalent:
    \begin{description}
      \item[(a)] System $[A; C]$ is ASMS;

      \item[(b)] For any $Q\in \mathscr{S}_n^{+}$, if  $(A,C|Q)$ is exactly observable (exactly detectable), then there exists a
      unique solution
      $S\in \mathscr{S}_n^{++}$ ($S\in \mathscr{S}_n^{+}$)  to the generalized Lyapunov equation (GLE)
          \begin{equation}\label{lyp111}
          A'SA+C'SC+Q=S;
          \end{equation}
      \item[(c)] The spectral set of ${\mathcal
D}_{A,C}$    satisfies    $\sigma ({\mathcal D}_{A,C})\subset$ $
{\mathcal D}(0,1):=\{\lambda: \lambda\in {\mathcal C},
|\lambda|<1\}$, where   the  generalized Lyapunov operator
${\mathcal D}_{A,C}$ is defined as
$$
{\mathcal  D}_{A,C}X=AXA'+CXC',    X\in {\mathscr S}_n.
$$
\end{description}
     \end{lemma}

The following PBH criteria can be found in
Theorem 3.7 of  \cite{zhang-xie-chen}.

\begin{lemma}[Stochastic PBH eigenvector
test]\label{the:Theorem3.3.2}

For  the  exact observability and exact detectability of  $(A, C|Q)$, we
have
\begin{itemize}
\item [(i)]  $(A, C|Q)$
is exactly observable if and only if (iff) there does not exist a non-zero $X\in
{\mathscr S}_n$  such that
\begin{equation}
{\mathcal D}_{A, C} X=\lambda X,\   CX=0, \  \lambda\in {\mathcal
C}. \label{eq the_3.3.2}
\end{equation}

\item [(ii)] $(A,C|Q)$
is exactly detectable iff  there does not exist a non-zero $X\in
{\mathscr S}_n$ such that
\begin{equation}
{\mathcal D}_{A,C} X=\lambda X,\ \  CX=0, \ \ |\lambda|\ge 1.
\label{eq the}
\end{equation}
\end{itemize}
\end{lemma}

Throughout this paper, we adopt the following  standard assumptions.
\begin{assumption}\label{assump1}
Assume that
\begin{description}
  \item[(1)] $Q\succeq 0$ and $R\succ0$;
  \item[(2)] System $[A, B; C, D]$ is stabilizable;
  \item[(3)] $(A,C|Q)$ is exactly observable or exactly detectable.
      \end{description}
\end{assumption}

For our convenient use, from now on, we consider the following cost
function
\begin{equation}\label{equgu}
\hat{J}(z_1\cdots,z_r, F)=\sum\limits_{l=1}^rJ(z_l, F)
\end{equation}
instead of  ${J}(z, F)$ as in  (\ref{cost2}), where $Z:=\sum
\limits_{l=1}^r z_lz_l'\succ 0$.  From the well-known LQ theory,
although  the optimal value of $J(z, F)$ depends on the
initial state $z$, the optimal feedback gain $F^*$ is unrelated to
$z$. Hence,
$$
F^*=\text{arg}\min\limits_{F\in\mathscr{F}}\hat{J}(z_1\cdots,z_r,
F)=\text{arg}\min\limits_{F\in\mathscr{F}}{J}(z, F).
$$

We state the concerned SLQR as follows:

{\bf Problem (SLQR)}. Solve the  following non-convex minimization
problem:
\begin{eqnarray}\label{slqr}
    \left\{\begin{array}{l}
{\hat J}(F^*):=\min\limits_{\substack{F}\in\mathscr{F}}\hat{J}(z_1\cdots,z_r, F)=\min\limits_{\substack{F}\in\mathscr{F}}\sum\limits_{l=1}^rJ(z_l, F), \\
\texttt{s.t.} \ \ x_{k+1}=(A+BF)x_k+(C+DF)x_kw_k,
    \end{array}
    \right.
\end{eqnarray}
where  $\mathscr{F}:=\{F: [A+BF; C+DF] \ \ \mbox {is ASMS}\}$.

Note that Assumption~\ref{assump1}-(2) guarantees that
the Problem (SLQR) is well-posed, while Assumption 2.1-(3) guarantees
that $F^*$  is also a feedback stabilizing gain, which makes the
closed-loop system ASMS \cite{huangyulin,zhang-xie-chen}.

From the work of \cite{huangyulin,zhang-xie-chen}, we have the
following results on  Problem (SLQR):

\begin{lemma}\label{lem2.X}
Under Assumption~\ref{assump1}, Problem (SLQR) is well-posed
and attainable, concretely speaking, the optimal gain $F^*$, which
minimizes the cost (\ref{equgu}), is given by
\begin{eqnarray}\label{Fdefine}
    F^* = -\left(R + B'P^*B + D'P^*D\right)^{-1}\left(B'P^*A + D'P^*C\right),
\end{eqnarray}
and the optimal value function  for the  Problem (SLQR) is
\begin{eqnarray}
    {\hat J}(F^*)=\sum_{l=1}^r z'_lP^* z_l=Tr(ZP^*),
\end{eqnarray}
where, under exact observability (exact detectability),   $P^*\in
{\mathscr S}^{++}$ $(P^*\in {\mathscr S}^{+})$ is the unique
solution to the GARE
\begin{eqnarray}\label{geare}
    A'PA + C'PC + Q - (A'PB + C'PD)\,(R + B'PB + D'PD)^{-1}\,(B'PA + D'PC) = P
\end{eqnarray}
\end{lemma}
Similar to (5) of \cite{lee2019},
 the Q-learning method provides a model-free solution
for solving SLQR.  Define Q-function  for SLQR  as
\begin{eqnarray}
    Q^*(x_k, u_k):&=&\mathcal{E}\{x_k'Qx_k+u'_kRu_k\}+\min\limits_{u}J(x_{k+1}, u)\nonumber\\
   & =&\mathcal{E}\left\{
    \left[\begin{array}{ccc}
    x_k\\u_k
    \end{array}
    \right]'X^*
    \left[\begin{array}{ccc}
    x_k\\u_k
    \end{array}
    \right]
    \right\},
\end{eqnarray}
where
\begin{eqnarray}
&&X^* \!=\!  \left[\begin{array}{ccccc}
   X^*_{11}&X^*_{12}\\
   (X^*_{12})'&X^*_{22}
    \end{array}
    \right]\nonumber\\
  &&  \ \ \  \!:=\!\left[\!\begin{array}{ccccc}
   Q\!+\!A'P^*A\!+\!C'P^*C\!&\!A'P^*B\!+\!C'P^*D\\
   B'P^*A\!+\!D'P^*C\!&\!R\!+\!B'P^*B\!+\!D'P^*D
    \end{array}\!
    \right].\nonumber\\
    \label{Xdefine}
\end{eqnarray}
The optimal control input  minimizes the Q-function, i.e.,
$$
u^*_k=F^*x_k=\mbox {argmin}_{u_k\in {\mathcal U}_{ad}}Q^*(x_k, u_k).
$$

\subsection{Off-line PI  for {Problem (SLQR)}}

This subsection
introduces an off-line PI to solve the  GARE
(\ref{geare}) arising from the  Problem (SLQR), which has  certain  connection
with the  primal-dual  model-free algorithm.

In order to give  a model-based PI method for  the  Problem
(SLQR), we give  the policy evaluation step and the policy update
step  as
\begin{eqnarray}\label{offlineP}
    P^{(i)} = (A+BF^{(i)})'P^{(i)}(A+BF^{(i)}) + (C+DF^{(i)})'P^{(i)}(C+DF^{(i)}) + Q + (F^{(i)})'RF^{(i)}
\end{eqnarray}
and
\begin{eqnarray}\label{offlineF}
 F^{(i+1)} = -\left(R + B'P^{(i)}B + D'P^{(i)}D\right)^{-1} \left(B'P^{(i)}A + D'P^{(i)}C\right)
\end{eqnarray}
respectively. It should be pointed out that the continuous-time
model-based PI  can be found in \cite{zhangphd}.

\begin{lemma}\label{converoff}
In the off-line PI algorithm, the two sequences
$\{P^{(i)}\}_{i=0}^\infty$ and $\{F^{(i)}\}_{i=0}^\infty$  have the
properties that
\begin{enumerate}
  \item $P^*\preceq P^{(i+1)}\preceq P^{(i)}$;
  \item   $\lim\limits_{i\rightarrow \infty}P^{(i)}=P^*$, $\lim\limits_{i\rightarrow \infty}F^{(i)}=F^*$,
  where $P^*$ is the solution to GARE (\ref{geare}) and $F^*$ is as given in (\ref{Fdefine}).
\end{enumerate}
\end{lemma}

{\bf Proof.} {\bf Step 1:} Because
 $[A,B; C,D]$  is stabilizable, there exists an $F^{0}$ such that  $[A+BF_0; C+DF_0]$ is
 ASMS. By  Theorem 3.7 of \cite{zhang-xie-chen},  it is easy to know that if $(A,C|Q)$ is exactly observable (exactly detectable),
 then so is $(A+BF^{(i)},C+DF^{(i)}|Q+(F^{(i)})'RF^{(i)})$.
 According to Lemma~\ref{lem2.1},  there exists a unique solution $P^{(0)}\in
\mathscr{S}_n^{++}$($P^{(0)}\in \mathscr{S}_n^{+}$) under exact
observability (exact detectability) for the policy evaluation
equation
 (\ref{offlineP}).

{\bf Step 2:}  In order to prove that $\{P^{(i)}\}$ can proceed for
ever and   is a monotonically decreasing sequence, as well as
$\{F^{(i)}\}$ is a feedback stabilizing gain sequence,  we need to
prove the following two iteration formulas: {\small
\begin{eqnarray}\label{iteration3}
    P^{(i)} = (A+BF^{(i+1)})'P^{(i)}(A+BF^{(i+1)}) + \tilde{Q} + (C+DF^{(i+1)})'P^{(i)}(C+DF^{(i+1)})
\end{eqnarray}
} and
\begin{align}\label{iteration4}
&(A+BF^{(i+1)})' \Delta P^{(i)}(A+BF^{(i+1)})
+(C+DF^{(i+1)})'\Delta P^{(i)}(C+DF^{(i+1)})-\Delta P^{(i)}\nonumber\\
=&-(\Delta F^{(i)})'(R+B'P^{(i)}B+D'P^{(i)}D)\Delta F^{(i)},
\end{align}
where
\begin{eqnarray*}\label{iteration5}
\tilde Q= Q+(F^{(i+1)})'RF^{(i+1)}+(\Delta F^{(i)})'(R+B'P^{(i)}B+D'P^{(i)}D)\Delta F^{(i)}
\end{eqnarray*} with $\Delta P^{(i)}=P^{(i)}-P^{(i+1)}$ and  $\Delta
F^{(i)}=F^{(i)}-F^{(i+1)}$.

Note that
\begin{align*}
&(A+BF^{(i+1)})' P^{(i)}(A+BF^{(i+1)})
+(C+DF^{(i+1)})' P^{(i)}(C+DF^{(i+1)})-P^{(i)}\\
=&(A+BF^{(i+1)})' P^{(i)}(A+BF^{(i+1)})
+(C+DF^{(i+1)})' P^{(i)}(C+DF^{(i+1)})
-(A+BF^{(i)})'P^{(i)}\\&\cdot(A+BF^{(i)})
-(C+DF^{(i)})'P^{(i)}(C+DF^{(i)})-Q-(F^{(i)})'RF^{(i)}\\
=&-(A'P^{(i)}B+C'P^{(i)}D)\Delta F^{(i)}
-(\Delta F^{(i)})'(B'P^{(i)}A+D'P^{(i)}C)-Q
+(F^{(i+1)})'(B'P^{(i)}B\\&+D'P^{(i)}D)F^{(i+1)}
-(F^{(i)})'(R+B'P^{(i)}B+D'P^{(i)}D)F^{(i)}\\
=&-(A'P^{(i)}B+C'P^{(i)}D)\Delta F^{(i)}
-(\Delta F^{(i)})'(B'P^{(i)}A+D'P^{(i)}C)-Q
+(F^{(i+1)})'\\&\cdot(R+B'P^{(i)}B+D'P^{(i)}D)F^{(i+1)}
-(F^{(i)})'(R+B'P^{(i)}B+D'P^{(i)}D)F^{(i)}
-(F^{(i+1)})'RF^{(i+1)}\\
=&F^{(i+1)}(R+B'P^{(i)}B+D'P^{(i)}D)\Delta F^{(i)}
+(\Delta F^{(i)})'(R+B'P^{(i)}B+D'P^{(i)}D) F^{(i+1)}-Q\\
&+(F^{(i+1)})'(R+B'P^{(i)}B+D'P^{(i)}D)F^{(i+1)}
-(F^{(i)})'(R+B'P^{(i)}B+D'P^{(i)}D)F^{(i)}\\
&-(F^{(i+1)})'RF^{(i+1)}\\
=&F^{(i+1)}(R+B'P^{(i)}B+D'P^{(i)}D)F^{(i)})
+F^{(i)})'(R+B'P^{(i)}B+D'P^{(i)}D) F^{(i+1)}-Q
-(F^{(i+1)})'\\&\cdot(R+B'P^{(i)}B+D'P^{(i)}D)F^{(i+1)}
-(F^{(i)})'(R+B'P^{(i)}B+D'P^{(i)}D)F^{(i)}
-(F^{(i+1)})'RF^{(i+1)}\\
=&-(\Delta F^{(i)})'(R+B'P^{(i)}B+D'P^{(i)}D)\Delta F^{(i)}-Q
-(F^{(i+1)})'RF^{(i+1)}.
\end{align*}
Hence, (\ref{iteration3}) is proved.  By combining \eqref{offlineP}
at step $i+1$ and \eqref{iteration3}, (\ref{iteration4}) is  easily
obtained.

Based on (\ref{offlineP}) and  (\ref{iteration3}), we know that if
$F^{(i)}$ is a feedback stabilizing gain matrix, then so is
$F^{(i+1)}$, which yields $P^{(i+1)}\succ 0$ ($P^{(i+1)}\succeq 0$)
under exact observability (exact detectability)  for $i\ge  0$
by repeating  Step 1. Hence, by solving
(\ref{offlineP})-(\ref{iteration3}), we can obtain a  sequence
$\{P^{(i)}\succ 0\}$ ($\{P^{(i)}\succeq 0\}$) under exact
observability (exact detectability)  and a feedback stabilizing gain
matrix sequence $\{F^{(i)}\}$ in the following order:
$$
F^{(0)}\rightarrow P^{(0)}\rightarrow F^{(1)}\rightarrow
P^{(1)}\rightarrow \cdots.
$$
By (\ref{iteration4}) and Lemma~\ref{lem2.1}, there exists a unique
solution $\Delta P^{(i)}\succ 0$ ($\Delta P^{(i)}\succeq  0$) to
(\ref{iteration4}), which results in that $\{P^{(i)}\succ 0\}$
($\{P^{(i)}\succeq 0\}$) is a monotonically decreasing sequence.

{\bf Step 3:} Because $\{P^{(i)}\succ 0\}$ ($\{P^{(i)}\succeq 0\}$)
is a monotonically decreasing sequence with low bound zero, it must
have a unique limit $P^*$ satisfying (\ref{geare}).   Moreover,
 $$
 \lim_{i\to\infty} P^{(i)}=P^*\preceq
P^{(i+1)}\preceq P^{(i)}.
$$
Taking the limit on both sides of (\ref{offlineF}), we have
\begin{eqnarray*}
\lim\limits_{i\rightarrow \infty}F^{(i+1)}
&=\lim\limits_{i\rightarrow\infty}[-(R+B'P^{(i)}B+D'P^{(i)}D)^{-1}(B'P^{(i)}A+D'P^{(i)}C)]\\
&=-(R+B'P^{*}B+D'P^{*}D)^{-1}(B'P^{*}A+D'P^{*}C)
:=F^*,
\end{eqnarray*}
where $F^*$ satisfies  (\ref{Fdefine}).
Combining  Step 1-Step 3, Lemma~\ref{converoff} is
proved. $\blacksquare$

\begin{remark}\label{eqvhvhaa}
Following the line of Theorem 2.4.1 of  \cite{zhangphd}, it is easy
to show that the convergence speed of $\{P^{(i)}\}_{i=0}^\infty$ is
quadratic, i.e.,
$$
\|P^{(i+1)}-P^{(i)}\|\le  C\|P^{(i)}-P^{(i-1)}\|^2,
$$
where  $C$ is a constant.
\end{remark}

\section{Primal-dual optimization-based method}
This section aims to explore the {Problem (SLQR)}
 from the perspective of primal-dual  optimization.
The primal-dual algorithms,  dependent  and independent of
the system model are investigated.

\subsection{Problem reformulation}

Considering technical requirements,  we construct the augmented state
vector as
$$
v_k:=\left[\begin{array}{ccc}x_k\\
u_k\end{array}\right], \ u_k=Fx_k.
$$
System (\ref{sys1}) with $x_0=z$  yields the following
augmented system:
\begin{eqnarray}\label{augv}
&&v_{k+1}=A_Fv_k+C_Fv_kw_k,\\
&&v_0=\left[\begin{array}{ccc}x_0\nonumber\\
Fx_0 \end{array}\right]=\left[\begin{array}{ccc}
I_n\\
F \end{array}\right]z\in \mathcal{R}^{n+m},
\end{eqnarray}
where
$$
A_F:=\left[\begin{array}{cccc}A&B\\FA&FB\end{array}\right]\in\mathcal{R}^{(n+m)\times(n+m)}
$$
and
$$
C_F:=\left[\begin{array}{cccc}C&D\\FC&FD\end{array}\right]\in
\mathcal{R}^{(n+m)\times(n+m)}.
$$
 The solution of the system (\ref{augv})
 concerning the  initial state $v_0$  is denoted as $v_k^{F,  v_0}$.
In particular,  if  $x_0=z$, $u_0=Fx_0=Fz$, then
 we also write $v_k^{F,  v_0}$ as $v_k^{F,  z}$  and
 $u_k=Fx_k^{F,z}$  as $u_k^{F,z}$.
Because
\begin{eqnarray*}
    \mathcal{E}\|v_k^{F, z}\|^2&=& \mathcal{E}\|x_k^{F,  z}\|^2+\mathcal{E}\|u_k^{F,  z}\|^2\\
    &=&\mathcal{E}\|x_k^{F,  z}\|^2+\mathcal{E}
    \|Fx_k^{F,  z}\|^2,
\end{eqnarray*}
which means that it is equivalent to that between the ASMS of the original
system (\ref{sysF}) and the ASMS of the  augmented system
(\ref{augv}).

\begin{lemma}\label{lem3.1}
System (\ref{sysF})  is  ASMS iff the  augmented system (\ref{augv})
is ASMS.
\end{lemma}

Under the constraint of the augmented system (\ref{augv}), the
associated
 cost function in (\ref{slqr})  can be  written as
 \begin{eqnarray}\label{eqvvh}
{\hat J}(z_1,\cdots, z_r,
F)=\mathcal{E}\sum\limits_{l=1}^r\sum\limits_{k=0}^\infty
(v_k^{F,z_l})'\Lambda (v_k^{F,z_l} ),
\end{eqnarray}
while the optimal  feedback  gain matrix  $F^*$ remains unchanged
\cite{zhang-xie-chen,huangyulin}.

In the following, the original  non-convex    optimization  problem
is transformed into  ($\mathcal{P}_1$-SLQR). The equivalence of
 ($\mathcal{P}_1$-SLQR)  to  Problem (SLQR) is also
provided.

{\bf   ($\mathcal{P}_1$-SLQR) (Primal Problem I).}
 Solve the following non-convex minimization problem:
\begin{eqnarray}
\left\{\begin{array}{l}
J_{\mathcal{P}_1}:=\inf\limits_{S\in\mathscr{S}_{n+m},\  F\in\mathcal{R}^{m\times n}} Tr(\Lambda S)\\
{\text{s.t.}} \
A_FSA_F'+C_FSC_F'+\left[\begin{array}{cccc}I_n\\F\end{array}\right]Z
\left[\begin{array}{cccc}I_n\\F\end{array}\right]'=S, \label{ly1}\\
F\in {\mathscr F}.\\
\end{array}\right.
\end{eqnarray}

\begin{proposition}
The optimization problem of { Problem (SLQR)} can be transformed
into that of { ($\mathcal{P}_1$-SLQR)}. Moreover,
$J_{\mathcal{P}_1}=J(F^*)$, $F_{\mathcal{P}_1}=F^*$.
\end{proposition}
{\bf Proof.}  Set
\begin{eqnarray*}\label{sdefi}
S&:=&\sum\limits_{l=1}^r\sum\limits_{k=0}^\infty \mathcal{E}[v_k^{F, z_l}(v_k^{F, z_l})'].
\end{eqnarray*}
Due to $F\in \mathscr{F}$, by Lemma~\ref{lem3.1},  we must have
$0\preceq S\prec \infty$, i.e., the above-defined $S$  is
meaningful. Then, the objective function in the Problem (SLQR) can
be equivalently expressed as
\begin{eqnarray}
\hat{J}(z_1,\cdots,z_r,
F)&=&\mathcal{E}\sum\limits_{l=1}^r\sum\limits_{k=0}^\infty
\left[\begin{array}{ccc}x_k^{F, z_l}\\Fx_k^{F,
z_l}\end{array}\right]'\Lambda
\left[\begin{array}{ccc}x_k^{F, z_l}\\Fx_k^{F, z_l}
\end{array}\right] \nonumber\\
&=&\mathcal{E}\sum\limits_{l=1}^r\sum\limits_{k=0}^\infty \left(v_k^{F,z_l}\right)'\Lambda v_k^{F, z_l}\nonumber\\
&=& Tr (  \Lambda S).
\end{eqnarray}
Along with system $[A_F; C_F]$, $k\geq 1$,  there is
\begin{eqnarray}\label{eqbjbn}
S
&\!=\!&\sum\limits_{l=1}^r\left\{\left[\begin{array}{ccc}I_n\\F\end{array}\right]z_lz_l'
\left[\begin{array}{ccc}I_n\\F\end{array}\right]'
\!+\!\sum\limits_{k=1}^\infty\mathcal{E}\left[v_k^{F, z_l}(v_k^{F, z_l})'\right]\right\}\nonumber\\
&\!=\!&\left[\begin{array}{ccc}I_n\\F\end{array}\right]Z
\left[\begin{array}{ccc}I_n\\F\end{array}\right]'
\!+\!\sum\limits_{l=1}^r\sum\limits_{k=1}^\infty \mathcal{V}_k^{F, z_l}
\end{eqnarray}
with $\mathcal{V}_k^{F, z_l}:=\mathcal{E}[v_k^{F, z_l}(v_k^{F,
z_l})']$. Since the vector $A_Fv_k$  is
$\mathcal{F}_{k-1}$-measurable and  is independent  of   $w_k$, one
obtains
\begin{eqnarray*}
&&\mathcal{E}\left[v_{k+1}^{F, z_l}(v_{k+1}^{F, z_l})'\right]
=A_F
\mathcal{E}\left[v_k^{F, z_l}(v_k^{F, z_l})']A_F'+C_F \mathcal{E}[v_k^{F, z_l}(v_k^{F, z_l})'\right]C_F'.
\end{eqnarray*}
Hence, from (\ref{eqbjbn}),
\begin{eqnarray*}
S&=&\left[\begin{array}{ccc}I_n\\F\end{array}\right]Z\left[\begin{array}{ccc}I_n\\F\end{array}\right]'
+A_F\left[\sum\limits_{l=1}^r\sum\limits_{k=0}^\infty
\mathcal{V}_k^{F,
z_l}\right]A_F'
+C_F\left[\sum\limits_{l=1}^r\sum\limits_{k=0}^\infty
\mathcal{V}_k^{F, z_l}\right]C_F'\\
&=&\left[\begin{array}{ccc}I_n\\F\end{array}\right]Z\left[\begin{array}{ccc}I_n\\F\end{array}\right]'+A_FSA_F'+C_FSC_F',
\end{eqnarray*}
that is,  $S$ satisfies the GLE (\ref{ly1}). From the
above discussion, it can be seen that
$J_{\mathcal{P}_1}=J(F^*)$ $=Tr(\Lambda S^*_{\mathcal{P}_1})$, where
$(S^*_{\mathcal{P}_1}, F^*)$ is
the unique solution  of (\ref{ly1})  with  $F^*$  being  given by (\ref{Fdefine}). $\blacksquare$

\begin{remark}
From Lemma \ref{lem2.1},   $F\in\mathscr{F} $ implies $[A+BF; C+DF]$
to be  ASMS. Hence, according to Lemma \ref{lem3.1},  $[A_F; C_F]$
is also  ASMS.
   By the result of \cite{Bouhtouri_1999},  (\ref{ly1})  admits  a  unique
   solution $S\succeq0$  for any $F\in\mathscr{F}$.
\end{remark}

If we express the  matrix $S$ in a block form as  $S:=\begin{bmatrix}S_{11}&S_{12}\\S_{12}'&S_{22}\end{bmatrix}$
   with $S_{11}\in\mathscr{S}_n$, $S_{12}\in\mathcal{R}^{n\times m}$, and $S_{22}\in\mathscr{S}_m$, the following properties hold.

\begin{lemma}
In {\bf   ($\mathcal{P}_1$-SLQR)},  any feasible solution
$S\in\mathscr{S}_{n+m}$ and   $F\in\mathcal{R}^{m\times n}$ satisfy
\begin{description}
  \item[1)] $(A+BF)S_{11}(A+BF)'+(C+DF)S_{11}(C+DF)'+Z=S_{11}$;
  \item[2)] $\begin{bmatrix}I_n\\F\end{bmatrix}S_{11}
  \begin{bmatrix}I_n\\F\end{bmatrix}'=S$;
    \item[3)] $F=S_{12}'S_{11}^{-1}$.
      \end{description}
\end{lemma}
 {\bf   Proof. }  Observe that the GLE
in (\ref{ly1})  can be rewritten as
\begin{eqnarray}\label{sdfsssss}
S&=&A_FSA_F'+C_FSC_F'+\left[\begin{array}{cccc}I_n\\F
\end{array}\right]Z
\left[\begin{array}{cccc}I_n\\F\end{array}\right]'\nonumber\\
& =&\left[\begin{array}{cccc}\!I_n\!\\\!F\!
\end{array}\right]\!\left[\begin{array}{cccc}\!A\!&\!B\!
\end{array}\right]\!S\!\left[\begin{array}{cccc}\!A'\!\\ \!B'\!
\end{array}\right]\!
\left[\begin{array}{cccc}\!I_n\!\\\!F\!
\end{array}\right]'\!+\!\left[\begin{array}{cccc}\!I_n\!\\ \!F\!
\end{array}\right]\!Z\!\left[\begin{array}{cccc}\!I_n\!\\ \!F\!
\end{array}\right]'\!\nonumber\\
&+& \left[\begin{array}{cccc}\!I_n\!\\\!F\!
\end{array}\right]\left[\begin{array}{cccc}\!C\!&\!D\!
\end{array}\right]S\left[\begin{array}{cccc}\!C'\!\\\!D'\!
\end{array}\right]
\left[\begin{array}{cccc}\!I_n\!\\\!F\!
\end{array}\right]'.
\end{eqnarray}
By comparing the first $n\times n$ block  matrix  of the above
equation, the following  holds:
\begin{equation}\label{s11eq}
S_{11}\!=\!\left[\begin{array}{cccc}\!A\!&\!B\!
\end{array}\right]\!S\!\left[\begin{array}{cccc}\!A'\! \\\!B' \!
\end{array}\right]\!+\! \left[\begin{array}{cccc}\!C\!&\!D\!
\end{array}\right]\!S\!\left[\begin{array}{cccc}\!C'\! \\\! D'\!
\end{array}\right]\!+\! Z.
\end{equation}
\allowdisplaybreaks
Using GLE  (\ref{ly1}) again, we have
\begin{eqnarray*}
\begin{array}{l}
S_{11}=\left[\begin{array}{cccc}\!A\!&\!B\!
\end{array}\right]\!
\left[\begin{array}{cccc}\!I_n\!\\\!F\!
\end{array}\right]\!
\left[\begin{array}{cccc}\!A\!&\!B\!
\end{array}\right]\!
S \!\left[\begin{array}{cccc}\!A'\!\\\!B'\!
\end{array}\right]\!
\left[\begin{array}{cccc}\!I_n\!&\!F'\!
\end{array}\right]\!
\left[\begin{array}{cccc}\!A'\!\\\!B'\!
\end{array}\right]\\
   \ \  +\left[\begin{array}{cccc}\!A\!&\!B\!
\end{array}\right]\!
\left[\begin{array}{cccc}\!I_n\!\\\!F\!
\end{array}\right]\!
\left[\begin{array}{cccc}\!C\!&\!D\!
\end{array}\right]\!
S \!\left[\begin{array}{cccc}\!C'\!\\ \!D'\!
\end{array}\right]\!
\left[\begin{array}{cccc}\!I_n\!&\!F'\!
\end{array}\right]\!
\left[\begin{array}{cccc}\!A'\!\\\!B'\!
\end{array}\right]
+\!\left[\begin{array}{cccc}\!A\!&\!B\!
\end{array}\right]
\left[\begin{array}{cccc}\!I_n\!\\\!F\!
\end{array}\right]\!Z\!
\left[\begin{array}{cccc}\!I_n\!&\!F'\!
\end{array}\right]\!
\left[\begin{array}{cccc}\!A'\!\\\!B'\!
\end{array}\right]\\
   \ \  +\!\left[\begin{array}{cccc}\!C\!&\!D\!
\end{array}\right]
\left[\begin{array}{cccc}\!I_n\!\\\!F\!
\end{array}\right]\!
\left[\begin{array}{cccc}\!A\!&\!B\!
\end{array}\right]\!
S\! \left[\begin{array}{cccc}\!A'\!\\\!B'\!
\end{array}\right]\!
\left[\begin{array}{cccc}\!I_n\!&\!F'\!
\end{array}\right]\!
\left[\begin{array}{cccc}\!C'\!\\\!D'\!
\end{array}\right]\\
  \ \  +\!\left[\begin{array}{cccc}\!C\!&\!D\!
\end{array}\right]\!
\left[\begin{array}{cccc}\!I_n\!\\\!F\!
\end{array}\right]\!
\left[\begin{array}{cccc}\!C\!&\!D\!
\end{array}\right]
S \left[\begin{array}{cccc}\!C'\!\\\!D'\!
\end{array}\right]\!
\left[\begin{array}{cccc}\!I_n\!&\!F'\!
\end{array}\right]\!
\left[\begin{array}{cccc}\!C'\!\\\!D'\!
\end{array}\right]\\
 \ \  + \!\left[\begin{array}{cccc}\!C\!&\!D\!
\end{array}\right]\!
\left[\begin{array}{cccc}\!I_n\!\\\!F\!
\end{array}\right]\!
Z \!\left[\begin{array}{cccc}\!I_n\!&\!F'\!
\end{array}\right]\!
\left[\begin{array}{cccc}\!C'\!\\\!D'\!
\end{array}\right]\!+\!Z\!\\
= \! (A\!+\!BF)\!\left[\begin{array}{cccc}\!A\!&\!B\!
\end{array}\right]\!S\!
\left[\begin{array}{cccc}\!A'\!\\\!B'\!
\end{array}\right]\!
(A\!+\!BF)'+\!(A\!+\!BF)\left[\begin{array}{cccc}\!C\!&\!D\!
\end{array}\right]\!S\!
\left[\begin{array}{cccc}\!C'\!\\\!D'\!
\end{array}\right]
(A\!+\!BF)'
\end{array}
\end{eqnarray*}
\begin{eqnarray*}
\begin{array}{l}
 \ \   +(A\!+\!BF)\!Z\!(A\!+\!BF)'
+\!(C\!+\!DF)\left[\begin{array}{cccc}\!A\!&\!B\!
\end{array}\right]\!S\!
\left[\begin{array}{cccc}\!A'\!\\\!B'\!
\end{array}\right]\!
(C\!+\!DF)'\\
 \ \   +\!(C\!+\!DF)\left[\begin{array}{cccc}\!C\!&\!D\!
\end{array}\right]\!S\!
\left[\begin{array}{cccc}\!C'\!\\\!D'\!
\end{array}\right]
(C\!+\!DF)'
+\!(C\!+\!DF)Z(C\!+\!DF)'\!+\!Z\\
=(A\!+\!BF)S_{11} (A\!+\!BF)'\!+\!(C\!+\!DF)S_{11} (C\!+\!DF)'+Z.
 \end{array}
\end{eqnarray*}
The first statement is obtained. Next, we plug the equivalence
relation for $S_{11}$ given in (\ref{s11eq}) into (\ref{sdfsssss}),
which leads to the second result. In addition,  the following
expanded form
\begin{eqnarray*}
S=\left[\begin{array}{cccc}I_n\\F
\end{array}\right]S_{11}\left[\begin{array}{cccc}I_n\\F
\end{array}\right]'
=\left[\begin{array}{cccc}S_{11}&S_{11}F'\\FS_{11}&FS_{11}F'
\end{array}\right]
\end{eqnarray*}
results in $S_{11}F'=S_{12}$. Since $S_{11} \succeq Z \succ  0$,
there is $F=S_{12}'S_{11}^{-1}$. The proof is completed.
$\blacksquare $

Below, based on our established    exact detectability and PBH
criteria, we are in a position to generalize Proposition 3 of
\cite{lee2019} to stochastic version.

\begin{proposition}\label{hhj}
The constrained condition $F\in {\mathscr F}$ in
($\mathcal{P}_1$-SLQR) can be replaced by $S\succeq 0$, where
$S\succeq 0$ is the unique solution of the GLE in  (\ref{ly1}).
\end{proposition}
{\bf    Proof. } By definition,  $F\in {\mathscr F}$ is equivalent to  that
$[A_F;C_F]$ is ASMS. By Lemma~\ref{lem2.1},  we only need to prove
that
$$
\left(A_F',C_F'|\tilde Q \right)
$$
is exactly detectable, where
$$
\tilde Q=\left[\begin{array}{cccc}I_n\\F\end{array}\right]Z
\left[\begin{array}{cccc}I_n\\F\end{array}\right]'.
$$
By Lemma~\ref{the:Theorem3.3.2}, we only need to show that there
does not have a non-zero $ X\in {\mathscr S}_{n+m}$ satisfying
\begin{equation}\label{eqvhvh}
A_F'XA_F+C_F'XC_F=\lambda X, \ |\lambda| \succeq 1
\end{equation}
and
\begin{equation}\label{eqvhvj}
\left[\begin{array}{cccc}I_n\\F\end{array}\right]Z
\left[\begin{array}{cccc}I_n\\F\end{array}\right]'X=0.
\end{equation}
Set $X=\begin{bmatrix}X_{11}&X_{ 12}\\X_{12}'&X_{22}
\end{bmatrix},
$ then (\ref{eqvhvj}) yields that
\begin{equation}\label{eqvhvaaj}
\begin{bmatrix}ZX_{11}+ZF'X_{12}'&ZX_{ 12}+ZF'X_{22}\\
FZX_{11}+FZF'X_{12}'&FZX_{12}+FZF'X_{22}
\end{bmatrix}=0.
\end{equation}
Because $Z\succ 0$, it can be derived from (\ref{eqvhvaaj}) that
\begin{equation}\label{eqgcf}
X_{11}+F'X_{12}'=0, \  X_{ 12}+F'X_{22}=0.
\end{equation}
In addition, considering (\ref{eqgcf}),  by computations, {\small
\begin{eqnarray}
&&A_F'XA_F=\begin{bmatrix}A'&A'F'\\
B'&B'F'
\end{bmatrix}\begin{bmatrix}X_{11}&X_{12}\\
X_{12}'&X_{22}
\end{bmatrix} \begin{bmatrix}A&B\\
FA&FB
\end{bmatrix}\nonumber\\
&&=\begin{bmatrix}A'&A'F'\\
B'&B'F'
\end{bmatrix}\begin{bmatrix}-F'X_{12}'&-F'X_{22}\\
-X_{22}F&X_{22}
\end{bmatrix} \begin{bmatrix}A&B\\
FA&FB
\end{bmatrix}\nonumber\\
&&=\begin{bmatrix}-A'F'X_{12}'-A'F'X_{22}F&-A'F'X_{22}+A'F'X_{22}\\
-B'F'X_{12}'-B'F'X_{22}F&-B'F'X_{22}+B'F'X_{22}
\end{bmatrix}\nonumber\\
&&=\begin{bmatrix} 0&0\\
0 &0
\end{bmatrix}.
\label{ewvhvh}
\end{eqnarray}
}
Similarly,
\begin{equation}\label{eqcvghv}
C_F'XC_F=\begin{bmatrix} 0&0\\
0 &0
\end{bmatrix}.
\end{equation}
Next, we show there does not have a non-zero $X\in {\mathscr
S}_{n+m}$ satisfying (\ref{eqvhvh}). In view of (\ref{ewvhvh}) and
(\ref{eqcvghv}), (\ref{eqvhvh}) holds iff
$$
\lambda \begin{bmatrix} X_{11}&X_{12}\\
X_{12}' & X_{22}
\end{bmatrix}=0, \ |\lambda| \succeq 1.
$$
{Without loss of generality, we assume
$\lambda=\lambda_1+i\lambda_2$
 with  $\lambda_1\ne 0$}. We only take $X_{11}$ as an example  to show $X_{ij}=0,
 i,j=1,2$. Assume $X_{11}=X_{11}^0+iX_{11}^1$. From $\lambda
 X_{11}=0$, it follows that
\begin{equation}\label{eqvhbj}
\lambda_1X_{11}^0-\lambda_2 X_{11}^1=0
\end{equation}
and
\begin{equation}\label{eqddvhbj}
\lambda_1X_{11}^1+\lambda_2 X_{11}^0=0.
\end{equation}
From  (\ref{eqvhbj}), $X_{11}^0=\frac
{\lambda_2}{\lambda_1}X_{11}^1$. Substitute the obtained $X_{11}^0$
into  (\ref{eqddvhbj}), we have
$$
\frac {\lambda_1^2+\lambda_2^2}{\lambda_1} X_{11}^1=0,
$$
which leads to $X_{11}^1=0$, and accordingly, $X_{11}^0=\frac
{\lambda_2}{\lambda_1}X_{11}^1=0$. So $X_{11}=0$. Repeating the same
procedures, we know $X_{12}=0$ and $X_{22}=0$. i.e., $X=0$. By
Lemma~\ref{the:Theorem3.3.2}-(ii),
$$
\left(A_F',C_F'|\tilde Q\right)
$$
is exactly detectable. By Lemma~\ref{lem2.1},  $F\in {\mathscr F}$,
i.e.,  $[A_F;C_F]$ is ASMS iff the GLE  in  (\ref{ly1}) admits a
unique solution $S\succeq 0$. $\blacksquare$

\begin{remark}\label{jbj}
From the proof of Proposition~\ref{hhj}, we can see that
$$
\left(A_F',C_F'|\tilde Q\right)
$$
is not exactly observable. This is because, when we take
$\lambda=0$, then any $X=(X_{ij})_{2\times 2}\ne0$ satisfying
\begin{equation}\label{eqgcf}
X_{11}+F'X_{12}'=0, \  X_{ 12}+F'X_{22}=0.
\end{equation}
 is a solution of the following equations.
\begin{equation}
{\mathcal D}_{A_F, C_F} X=\lambda X,\   \tilde Q X=0.
\end{equation}

\end{remark}

For further analysis, the following optimal control  problem is
introduced, which can help
 to obtain a strong duality theorem.

 {\bf   ($\mathcal{P}_2$-SLQR) (Primal Problem II).}
 Solve the non-convex optimization  with variables
$X\in\mathscr{S}_{n+m}$, $F\in\mathcal{R}^{m\times n}$:
\begin{eqnarray}
\left\{\begin{array}{l}
 {J}_{\mathcal{P}_2}:=\inf\limits_{X\in {\mathscr S}_{n+m}, \ F\in {\mathcal R}^{m\times n}} Tr\left(\begin{bmatrix}I_n\\F\end{bmatrix}Z\begin{bmatrix}I_n\\F\end{bmatrix}'X\right),\\
 \text{s.t.} \
A_F'XA_F+C_F'XC_F+\Lambda-X=0,\\
F\in {\mathscr F}. \label{lyp2}
\end{array}\right.
\end{eqnarray}
The following proposition shows the equivalence of the initial
{Problem (SLQR)} and
{($\mathcal{P}_2$-SLQR)}.
\begin{proposition}\label{pro32}
The ($\mathcal{P}_2$-SLQR) has a unique optimal solution
$(X^*_{\mathcal{P}_2}$,  $F^*_{\mathcal{P}_2})$,  and it  is equivalent
to {Problem (SLQR)}   in the sense that
${J}_{\mathcal{P}_2}=J(F^*)$ $={J}_{\mathcal{P}_1}$ and
$F^*_{\mathcal{P}_2}=F^*$.
\end{proposition}
{\bf Proof.} By Lemma~\ref{lem3.1} and  $F\in \mathscr{F}$,
$[A_F;C_F]$ is ASMS. Hence, for any  $F\in \mathscr{F}$,  the GLE in
(\ref{lyp2}) admits a unique solution $X\in\mathscr{S}^+_{n+m}$ due to
$\Lambda\ge 0$. It is easy to verify that for any  $F\in
\mathscr{F}$, we have
\begin{eqnarray*}
&&{\hat J}(z_1,z_2, \cdots,z_r,F)\\
&=&\sum\limits_{l=1}^r \sum\limits_{k=0}^\infty \mathcal{E}
\left[({v}_{k}^{F, z_l})' \Lambda
{v}_{k}^{F, z_l}\right]\\
&=&\sum\limits_{l=1}^r \sum\limits_{k=0}^\infty \mathcal{E}
\left[({v}_{k}^{F, z_l})' \Lambda {v}_{k}^{F, z_l}
+\Delta V_{k}^{F, z_l}\right]
+\sum\limits_{l=1}^r\left\{ \mathcal{E} [({v}_{0}^{F,z_l})' X
{v}_{0}^{F,z_l}]-\lim\limits_{k\rightarrow \infty}
\mathcal{E}[({v}_{k}^{F,z_l})' X {v}_{k}^{F,z_l}]\right\}\\
&=&\sum\limits_{l=1}^r \sum\limits_{k=0}^\infty \mathcal{E}
\left[({v}_{k}^{F, z_l})' \prod
{v}_{k}^{F, z_l}\right] +\sum\limits_{l=1}^r \mathcal{E} [({v}_{0}^{F,z_l})' X
{v}_{0}^{F,z_l}]\\
&=&\sum\limits_{l=1}^r \mathcal{E} [({v}_{0}^{F,z_l})' X
{v}_{0}^{F,z_l}]
\end{eqnarray*}
where  $\prod:=A_F'XA_F+C_F'XC_F+\Lambda-X=0$, $\Delta
V_{k}^{F,z_l}=
\mathcal{E}\left[({v}_{k+1}^{F,z_l})'X{v}_{F,k+1}^{z_l}\right]-\mathcal{E}\left[({v}_{k}^{F,z_l})'
X {v}_{k}^{F,z_l}\right]$. Considering $X$ and $F$  satisfy    GLE in
(\ref{lyp2}), we have
\begin{eqnarray*}
{\hat J}(z_1,\cdots,z_r,F^*)&=& \inf\limits_{X\in {\mathscr S}_{n+m}, F\in
\mathscr{F}}\sum\limits_{l=1}^r
\mathcal{E} \left[({v}_{0}^{F, z_l})' X {v}_{0}^{F,z_l}\right]\\
&=&\inf\limits_{X\in {\mathscr S}^+_{n+m}, F\in \mathscr{F}}
Tr\left[\begin{bmatrix}I_n\\F\end{bmatrix}Z\begin{bmatrix}
I_n\\F\end{bmatrix}'X\right]\\
&=&{J}_{\mathcal{P}_2}.
\end{eqnarray*}
Similar to Proposition~\ref{hhj}, by Lemma~\ref{the:Theorem3.3.2},  if $(A,C|Q)$ is exactly observable (exactly detectable), then so is $(A_F,C_F|\Lambda)$. By Lemma~\ref{lem2.1} and the well-known SLQR theory \cite{zhang-xie-chen}, GLE in
(\ref{lyp2}) admits a unique solution $(X^*_{\mathcal{P}_2},F^*_{\mathcal{P}_2})$, which is the unique optimal solution of ($\mathcal{P}_2$-SLQR).
Because $F^*\in {\mathscr F}$ also makes ${J}_{\mathcal{P}_2}$ arrive at
the minimum, we must have $F^*_{\mathcal{P}_2}=F^*$.  $\blacksquare$

It can be easily shown that the solution $X$  of GLE  in
(\ref{lyp2}) can be expressed as $X=\sum\limits_{k=0}^\infty Y_k, $
 where
$Y_k$  solves
\begin{eqnarray*}
\left\{\begin{array}{l}
Y_{k+1}=A_F'Y_kA_F+C_F'Y_kC_F,\\
Y_0=\Lambda.
\end{array}\right.
\end{eqnarray*}

Denote
$
X^*_{\mathcal{P}_2}=\begin{bmatrix}X^*_{\mathcal{P}_2,  11}&X^*_{\mathcal{P}_2, 12}\\(X^*_{\mathcal{P}_2, 12})'&X^*_{\mathcal{P}_2, 22}
\end{bmatrix}
$ as the unique solution of
the GLE in  (\ref{lyp2}), the following property can be
directly obtained based on (\ref{Fdefine})  and  (\ref{Xdefine}),
and we omit the proof for simplicity.

\begin{lemma}\label{lem33}
The optimal solution $(X^*_{\mathcal{P}_2}, F^*_{\mathcal{P}_2})$ for {\bf   ($\mathcal{P}_2$-SLQR)} satisfies
\begin{description}
  \item[(a)] $X^*_{\mathcal{P}_2}=X^*$;
  \item[(b)] $X^*_{\mathcal{P}_2}\!\succeq \!0, X^*_{\mathcal{P}_2, 22}\!\succ\! 0,
  F^*_{\mathcal{P}_2}\!=\!-(X^*_{\mathcal{P}_2, 22})^{-1} \!(X^*_{\mathcal{P}_2, 12})'$.
      \end{description}
\end{lemma}

Associated with  ($\mathcal{P}_1$-SLQR),  the dual problem  is defined
as follows:

{\bf Problem ($\mathcal{D}$-SLQR) (Dual Problem).} Solve
\begin{eqnarray*}
J_{\mathcal{D}}&:=&\sup\limits_{X\in\mathscr{S}_{n+m}}d(X)\\
& =& \sup\limits_{X\in\mathscr{S}_{n+m}, }
\inf\limits_{{S}\in\mathscr{S}_{n+m}^{+}, F\in\mathscr{F}}L(X, F,
{S}),
\end{eqnarray*}
where
\begin{eqnarray*}
&& \ \ \ \L(X, F, {S})\\
&&=Tr(\Lambda {S})
+Tr\Big[\Big(\!A_F{S}A_F'\!+\!C_F{S}C_F'\!-\!
{S}\!+\!\left[\begin{array}{cccc}\!I_n\!\\\!F\!\end{array}\right]\!Z\!
\left[\begin{array}{cccc}\!I_n\!\\\!F\!\end{array}\right]'\Big)\!X\!\Big]\\
&& = Tr
\left(Z\left[\begin{array}{cccc}I_n\\F\end{array}\right]'X\left[\begin{array}{cccc}I_n\\F\end{array}\right]\right)
+Tr\left((A_F'{X}A_F+C_F'{P}C_F-X+\Lambda)S\right),
\end{eqnarray*}
and $d(X):=\inf\limits_{{S}\in\mathscr{S}_{n+m}^{+},
F\in\mathscr{F}}L(X, F, {S})$ is  the Lagrangian dual function. By
weak  duality theorem  \cite{boyd2004}, $J_{\mathcal D}\preceq
J_{{\mathcal P}_1}$. Below, we prove that the strong duality
still holds.

\begin{theorem}\label{strongd}
(Strong Duality) The following equality  holds:
$$
{J}_{\mathcal{P}_1}=J_{\mathcal{D}}.
$$
\end{theorem}
{\bf Proof.} Similar to Theorem 1 of \cite{lee2019}, set $\mathcal{X}: =\{X\in\mathscr{S}^{+}_{n+m}:
A_F'XA_F+C_F'XC_F-X+\Lambda \succeq 0, \forall F\in\mathscr{F}\}$,
then the Lagrangian dual function is given as
\begin{eqnarray*}
d(X)&=&\inf\limits_{{F\in \mathscr{F}, \  {S}\in\mathscr{S}_{n+m}^{+}}}L(X, F, {S})\\
&=&\left\{\begin{array}{l} \inf\limits_{F\in \mathscr {F}}
Tr\left(Z\left[\begin{array}{cccc}I_n\\F\end{array}\right]'X
\left[\begin{array}{cccc}I_n\\F\end{array}\right]\right), \ X\in\mathcal{X},\\
-\infty,\ \ \ \ \ otherwise.
\end{array}\right.
\end{eqnarray*}
 Next, we need to prove that ${\mathcal{X}}$ is nonempty. In fact, the optimal solution   $X^*_{\mathcal{P}_2}$ of
 ($\mathcal{P}_2$-SLQR)  belongs to ${\mathcal{X}}$.
By  Lemma \ref{lem33}, there is
\begin{eqnarray*}
A_{F^*_{\mathcal{P}_2}}'X^*_{\mathcal{P}_2}A_{F^*_{\mathcal{P}_2}}
+C_{F^*_{\mathcal{P}_2}}'X^*_{\mathcal{P}_2}
C_{F^*_{\mathcal{P}_2}}+\Lambda=X^*_{\mathcal{P}_2},
\end{eqnarray*}
where
\begin{eqnarray*}
F^*_{\mathcal{P}_2}=-(X^*_{\mathcal{P}_2, 22})^{-1} (X^*_{\mathcal{P}_2,
12})'.
\end{eqnarray*}
Then, by Lemma \ref{lem33}-(b) and Lemma 4 in \cite{lee2019}, there
is
\begin{eqnarray*}
\begin{bmatrix}I_n\\F\end{bmatrix}'X^*_{\mathcal{P}_2}\begin{bmatrix}I_n\\F
\end{bmatrix}&\succeq& X^*_{\mathcal{P}_2, 11}-X^*_{\mathcal{P}_2, 12}
(X^*_{\mathcal{P}_2, 22})^{-1}(X^*_{\mathcal{P}_2, 12})'\\
&=&\begin{bmatrix}I_n\\F^*_{\mathcal{P}_2}\end{bmatrix}'
X^*_{\mathcal{P}_2}
\begin{bmatrix}I_n\\ F^*_{\mathcal{P}_2}\end{bmatrix}.
\end{eqnarray*}
Because
$$
A_F'X^*_{\mathcal{P}_2}A_F=\begin{bmatrix}A & B\end{bmatrix}'\begin{bmatrix}I_n\\F\end{bmatrix}'X^*_{\mathcal{P}_2}\begin{bmatrix}I_n\\F
\end{bmatrix}\begin{bmatrix}A & B\end{bmatrix}
$$
and
$$
C_F'X^*_{\mathcal{P}_2}C_F=\begin{bmatrix} C & D \end{bmatrix}'\begin{bmatrix}I_n\\F\end{bmatrix}'X^*_{\mathcal{P}_2}\begin{bmatrix}I_n\\F
\end{bmatrix}\begin{bmatrix} C & D\end{bmatrix},
$$
we have
$A_F'X^*_{\mathcal{P}_2}A_F+C_F'X^*_{\mathcal{P}_2}C_F +\Lambda\succeq
A_{F^*_{\mathcal{P}_2}} X^*_{\mathcal{P}_2}A_{F^*_{\mathcal{P}_2}}
+C_{F^*_{\mathcal{P}_2}}X^*_{\mathcal{P}_2}C_{F^*_{\mathcal{P}_2}}
+\Lambda=X^*_{\mathcal{P}_2}$ for all $F\in\mathscr{F}$, which means
that $X^*_{\mathcal{P}_2} \in {\mathcal{X}}$. Therefore, the dual
problem is equivalent to
\begin{eqnarray*}
J_{\mathcal{D}}&=&\sup\limits_{X\in\mathscr{S}_{n+m}}
d(X)\\
&=&\sup\limits_{X\in{\mathcal{X}}}\inf_{F\in {\mathscr F}}
Tr\left(Z\begin{bmatrix}I_n\\F\end{bmatrix}'X
\begin{bmatrix}I_n\\F\end{bmatrix}\right).
\end{eqnarray*}
For $X_{\mathcal{P}_2} \in {\mathcal{X}}$, it follows that
\begin{equation}\label{eqhuh}
d(X_{\mathcal{P}_2})=\inf_{F\in {\mathscr F}}
Tr\left(Z\begin{bmatrix}I_n\\F \end{bmatrix}'X_{\mathcal{P}_2}
\begin{bmatrix}I_n\\F\end{bmatrix}\right).
\end{equation}
By definition, $d(X_{\mathcal{P}_2})\preceq J_{\mathcal{D}}$. Since
$X_{\mathcal{P}_2}\succeq 0$, and the objective function of
(\ref{eqhuh}) is quadratic with respect to $F$, the infimum of
(\ref{eqhuh}) is arrived at
$$
F^*_{\mathcal{P}_2}=-(X^*_{\mathcal{P}_2, 22})^{-1} (X^*_{\mathcal{P}_2,
12})'.
$$
As discussed in \cite{lee2019},  we must have
\begin{equation}\label{eqhvh}
{J}_{\mathcal{P}_2}=d(X^*_{\mathcal{P}_2})\preceq
J_{\mathcal{D}}\preceq {J}_{\mathcal{P}_1}.
\end{equation}
By Proposition~\ref{pro32},
${J}_{\mathcal{P}_2}={J}_{\mathcal{P}_1}$. Hence, (\ref{eqhvh})
yields $J_{\mathcal{D}}={J}_{\mathcal{P}_1}$, and the strong duality
 holds.  $\blacksquare$

In order to find a positive definite solution to (\ref{ly1}), a modification of
($\mathcal{P}_1$-SLQR) is made by defining
\begin{equation}\label{sdefi2}
\widetilde{S}=\sum\limits_{l=1}^r \sum\limits_{k=0}^\infty
\mathcal{E}[v_k^{F,  v_0^l}(v_k^{F, v_0^l})']
\end{equation}
with the initial state satisfying
\begin{equation}\label{initialkkk}
z_lz_l'=\begin{bmatrix}\!A\!& \!B
\end{bmatrix}
v_0^l(v_0^l)'\begin{bmatrix}\!A\!&\! B
\end{bmatrix}'+\begin{bmatrix}\!C\!& \!D\!
\end{bmatrix}v_0^l(v_0^l)'\begin{bmatrix}\!C\!&\! D\!
\end{bmatrix}'
\end{equation}
and $\Xi=\sum\limits_{l=1}^r v_0^l(v_0^l)'\succ 0$  for some $r\in
\{1,2,\cdots,\}$. Here $x_0=z_1,z_2,\cdots, z_r$, $u_0=u_0^1$,
$u_0^2$, $\cdots$, $u_0^r$, $ v_0^l=\begin{bmatrix}z_l\\ u_0^1
\end{bmatrix}.
$
Note that the initial control input  in  $u_0^l$ is freely chosen,
in other words,  the feedback form only starts  from  $k=1$.
Therefore, $\Xi=\sum\limits_{l=1}^r v_0^l(v_0^l)'\succ 0$ can always
be  achieved by choosing suitable $z_l$ and $u_0^l$.

Define
\begin{eqnarray*}
{\hat J}(v^l_0,F):=\mathcal{E}\left\{\sum_{k=0}^\infty (v_k^{F, v^l_0})'
    \Lambda v_k^{F, v^l_0} \right\}
\end{eqnarray*}
and
\begin{eqnarray*}
{J}(z_l,F):=\mathcal{E}\left\{\sum_{k=0}^\infty (v_k^{F, z_l})'
    \Lambda v_k^{F,z_l} \right\}
\end{eqnarray*}
where
$$
v_0^{F,v_0^l}=v^l_0=\begin{bmatrix}x^l_0\\
 u^l_0
\end{bmatrix}=\begin{bmatrix}z_l\\
u^l_0
\end{bmatrix}.
$$
and
$$
v_k^{F,v_0^l}=\begin{bmatrix}x_k\\
Fx_k
\end{bmatrix}=\begin{bmatrix}I_n\\
F
\end{bmatrix}x_k, \  k\in\{1,2,\cdots\}.
$$
It is easy to show that, if the initial value satisfies
(\ref{initialkkk}), we have
$$
{\hat J}(v_0^l,F)={J}(z_l,F)+(v_0^l)'\Lambda v^l_0.
$$
Hence,
$$
\sum_{l=1}^r {\hat J}(v_0^l,F)=\sum_{l=1}^r {J}(z_l,F)+\sum_{l=1}^r
(v_0^l)'\Lambda v_0^l.
$$

{\bf   ($\hat{\mathcal{P}}_1$-SLQR) (Modified Primal Problem I).}
 Solve the following non-convex minimization with variables $\widetilde{S}\in\mathscr{S}_{n+m}$ and  $F\in\mathcal{R}^{m\times
 n}$:
\begin{eqnarray}
\left\{\begin{array}{l}
\hat{J}_{\mathcal{P}_1}:=\inf\limits_{\widetilde{S}\in {\mathscr S}_{n+m}, \ F\in {\mathcal R}^{m\times n}}
 Tr(\Lambda \widetilde{S}),  \\
{\text{s.t.}} \
 A_F\widetilde{S}A_F'+C_F\widetilde{S}C_F'
 +\Xi=\widetilde{S}, \label{lymode}\\
\ \ \ \ \  \widetilde{S} \succ 0
\end{array}\right.
\end{eqnarray}
with $\Xi=\sum\limits_{l=1}^r v_0^l(v_0^l)'\succ 0$.

Note that   (${{\mathcal{P}}_2}$-SLQR) should be changed as
\begin{eqnarray}
\left\{\begin{array}{l}
\hat {J}_{{\mathcal{P}}_2}=\inf\limits_{{{X}\in {\mathscr S}_{n+m}, \ F\in {\mathcal R}^{m\times n}}}Tr\Big( \Xi X\Big),  \\
{\text{s.t.}} \
 A'_FXA_F+C'_FXC_F
 +\Lambda-X=0, \label{lymode}\\
\ \ \ \ \  X\succeq 0.
\end{array}\right.
\end{eqnarray}
It is similarly shown that
$\hat{J}_{{\mathcal{P}}_1}=\hat{J}_{{\mathcal{P}}_2}$. From
Proposition 1 and Proposition 3 in \cite{liman}, we can directly
know that the optimal solution $(X^*_{{\mathcal{P}}_2},
F^*_{{\mathcal{P}}_2})=({\hat X}^*_{{\mathcal{P}}_2},  {\hat
F}^*_{{\mathcal{P}}_2})$.

\begin{remark}
 Based on the same discussion as Proposition 6 in \cite{lee2019} and
 Proposition 1 in \cite{liman}, no matter the initial state is
$x_0=z$ or $v_0$  with   $u_0$ being freely
chosen,  the optimal feedback gain $F^*$ remains  unchanged.
\end{remark}

\begin{remark}
Because for SLQR problem, the optimal feedback gain $F^*$ does not
depend on the initial state and the initial time, if we only aim to
search for $F^*$,  the initial constrained condition
(\ref{initialkkk}) can be ignored.
\end{remark}

For the ($\hat{\mathcal{P}}_1$-SLQR), the Lagrangian function takes
the form of
\begin{eqnarray*}
{\hat L}(X, F, \widetilde{S}, X_0)=Tr\Big(X \Xi\Big)+Tr\Big[\Big(A_F'XA_F+C_F'XC_F-X-X_0+\Lambda)\Big)\widetilde{S}\Big]
\end{eqnarray*}
with $X\in {\mathscr S}_{n+m}$ and $X_0\in {\mathscr S}^+_{n+m}$,
and the dual problem is
\begin{eqnarray*}
{\hat J}_{\mathcal{D}}\!=\!\sup\limits_{X\in\mathscr{S}_{n+m},
X_0\in\mathscr{S}^+_{n+m}}
\inf\limits_{\widetilde{S}\in\mathscr{S}_{n+m}^{++}, F\in
{\mathcal{R}}^{m\times n}}{\hat L} (X, F, \widetilde{S},X_0).
\end{eqnarray*}
Similar to Theorem~\ref{strongd}, we have the following modified
strong dual theorem.

\begin{theorem}\label{strongddd}
(Strong Duality) The strong duality  for ($\hat{\mathcal{P}}_1$-SLQR)
holds, that is,
$$
\hat{J}_{\mathcal{P}_1}={\hat J}_{\mathcal{D}}.
$$
\end{theorem}

Below,  we derive the KKT conditions for primal problem  ($\hat{\mathcal{P}}_1$-SLQR) based on strong dual theorem, which plays a critical role in the  partially model-free design.

\begin{proposition}\label{prokkt}
If $(\widetilde{S}^*, F^*)$ and $(X^*, X_0^*)$  are
 the    primal  and dual optimal points of  ($\hat{\mathcal{P}}_1$-SLQR), respectively,
 then
 $(\widetilde{S}^*, F^*, X^*, X_0^*)$ satisfies the following KKT conditions:
\begin{eqnarray}
&&A_F\widetilde{S}A_F'+C_F\widetilde{S}C_F'+\Xi-\widetilde{S}=0,\label{kkt1}\\
&&\widetilde{S}\succ 0, \label{kkt2}\\
&&A_F'XA_F+C_F'XC_F+\Lambda-X=0, \label{kkt3}\\
&&(X_{12}'\!+\!X_{22}F)\Big(\!\begin{bmatrix}\!A\!&\!B\!\end{bmatrix}
\!\widetilde{S}\!\begin{bmatrix}\!A\!&\!B\!\end{bmatrix}'\!+\!\begin{bmatrix}
\!C\!&\!D\!\end{bmatrix}\!\widetilde{S}\!\begin{bmatrix}\!C\!&\!D\!
\end{bmatrix}'\!\Big)=0.\label{kkt4}
\end{eqnarray}
\end{proposition}
{\bf Proof.} By Theorem~\ref{strongddd}, the strong dual theorem
holds. Hence, KKT condition is satisfied, i.e., we have the
following:

{Feasible condition}:
\begin{eqnarray}
&&
A_{F^*} \widetilde{S}^*(A_{F^*})'+C_{F^*}\widetilde{S}^*(C_{F^*})'+
\Xi-\widetilde{S}^*=0,\\
\label{ed}
&&
{\widetilde S}^*\succ0,\\
\label{edd}
&&
X^*_0\succeq 0.
\label{eddd}
\end{eqnarray}
{ Complementary slackness condition:}
\begin{equation}
Tr(\widetilde{S}^*X^*_0)=0.
\label{edddd}
\end{equation}
{The stationary condition:}

\begin{eqnarray}
&&\left.\frac{\partial {\hat L}(X^*, F^*, \widetilde{S}, X_0^*)}{\partial \widetilde{S}}\right|_{\widetilde{S}=\widetilde{S}^*}\nonumber\\
&&=(A_F^*)'X^*A_{F^*}+(C_{F^*})'X^*C_{F^*}-X^*-X^*_0+\Lambda=0,\\
\label{stationary1}
&&\frac{\partial {\hat L}(X^*, F, \widetilde{S}^*, X_0^*)}{\partial F}\Big|_{F=F*}\nonumber\\
&&=
\left.\frac{\partial  Tr\Big\{\Big(\begin{bmatrix}\!I_n\!&\! F'\!\end{bmatrix}'\!\Psi\!
\begin{bmatrix}\!I_n\!& \!F'\!\end{bmatrix}\Big)\!X^*\!\Big\}
}{\partial F}\right|_{F\!=\!F*}\!=\!0\label{stationary2}
\end{eqnarray}
with
$$
\Psi =
 \begin{bmatrix}A&B\end{bmatrix}\widetilde{S}^* \begin{bmatrix}A&B\end{bmatrix}'+\begin{bmatrix}C&D
 \end{bmatrix}\widetilde{S}^* \begin{bmatrix}C&D\end{bmatrix}'.
$$
Since
\begin{eqnarray*}
&&Tr\Big\{\Big(\begin{bmatrix}I_n&F'\end{bmatrix}'\Psi
\begin{bmatrix}I_n&F'\end{bmatrix}\Big)X^*\Big\}\\
&&=Tr\Big(\begin{bmatrix}\Psi&\Psi F'\\ F\Psi&F\Psi F'\end{bmatrix}
\begin{bmatrix}X^*_{11}&X^*_{12}\\ (X^*_{12})'&X^*_{22}\end{bmatrix}\Big)\\
&&=Tr\Big(\Psi X^*_{11}+\Psi F' (X^*_{12})'+F\Psi X^*_{12}+F'\Psi
F'X^*_{22}\Big),
\end{eqnarray*}
(\ref{stationary2}) leads to
{\small
\begin{eqnarray}
\left.\frac{\partial {\hat L}(X^*, F, \widetilde{S}^*, X_0^*)}{\partial
F}\right|_{F=F^*}=2[(X^*_{12})'+X_{22}F^*]\Psi=0,
\label{station}
\end{eqnarray}
}
and  (\ref{kkt4}) is derived  accordingly. Combing (\ref{edd}), (\ref{eddd}) and (\ref{edddd}), we know $X_0^*=0$.
In view of (\ref{ed}), (\ref{edd}), (\ref{stationary1}) and (\ref{station}), (\ref{kkt1})-(\ref{kkt4}) are obtained. This proposition is shown. $\blacksquare$

\begin{remark}
By Lemma~\ref{lem33}-(b), (\ref{kkt4}) can be replaced by
\begin{equation}
F^*=-(X^*_{22})^{-1}(X^*_{12})'. \label{eqjj}
\end{equation}
In the following subsection, we will find that, the KKT condition (
(\ref{kkt1})-(\ref{kkt3}) and (\ref{eqjj})) makes it possible to
construct model-based and partially model-free  SLQR control policy
design.

\end{remark}

\subsection{Primal-dual algorithm}

In this subsection, both the model-based and the partially
model-free primal-dual algorithms are  introduced  for the SLQR
design problem.
 The   convergence analysis of the algorithm  reveals the  connections among the primal-dual  algorithm,
the  classical  PI  and  Q-leaning algorithm.

The model-based  procedure  for solving the KKT condition is  given
in Algorithm 1. In particular,  $(X^i, F^i)$ in Algorithm 1
converges to the optimal value
 $(X^*, F^*)$  defined in (\ref{Xdefine}) and (\ref{Fdefine}).

 \begin{algorithm}
  \caption{Model-Based Primal-Dual Algorithm}
  \begin{algorithmic}[1]  
          \State {Initialization:} $F^0\in\mathscr{F}$,  the convergence tolerance
          $\varepsilon >0$, and the initial iteration  $i=0$;
           \State  {Repeat;}
            \State  {Dual update:} Solve $X^i$ from the equation
            \begin{equation}\label{dualup}
            (A_{F^i})'X^iA_{F^i}+(C_{F^i})'X^iC_{F^i}+\Lambda=X^i;
            \end{equation}
             \State   {Primal update:}
               \begin{equation}\label{primalupw}
             F^{i+1}=-(X_{22}^i)^{-1}(X_{12}^i)';
             \end{equation}
             \State  $i\leftarrow i+1$;
             \State {Until} $\|F^i-F^{i+1}\|\leq \varepsilon$.
    \end{algorithmic}
\end{algorithm}

\begin{theorem}\label{thddddddd}
For  the two sequences $\{X^i\}_{i=0}^\infty$ and
$\{F^i\}_{i=0}^\infty$ in   Algorithm 1,
 there are $\lim\limits_{i\rightarrow \infty} X^i=X^*$  and
$\lim\limits_{i\rightarrow \infty}F^i=F^*$, where $X^*$ and $F^*$ are defined in (\ref{Xdefine}) and (\ref{Fdefine}), respectively.
\end{theorem}
{\bf Proof.}  Notice that  pre-and post-multiplying (\ref{dualup})
by $\begin{bmatrix}I_n&(F^i)'\end{bmatrix}$ and its transpose, there
is
\begin{eqnarray*}
P^i&&=(A+BF^i)'P^i(A+BF^i)+Q+(F^i)'RF^i
+(C+DF^i)'P^i(C+DF^i)\\
&&=(F^i)'(B'P^iB+D'P^iD+R)F^i+(F^i)'(B'P^iA
+D'P^iC)+(A'P^iB+C'P^iD)F^i\\
&&+(A'P^iA+C'P^iC+Q)
\end{eqnarray*}
with
$$
P^i:=\begin{bmatrix}I_n&(F^i)'\end{bmatrix}X^i\begin{bmatrix}I_n&(F^i)'\end{bmatrix}',
$$
which is equivalent to the policy evaluation step (\ref{offlineP})
in off-line PI algorithm. We further expand
$$
\begin{bmatrix}I_n&(F^i)'\end{bmatrix}X^i\begin{bmatrix}I_n&(F^i)'\end{bmatrix}'
$$
with
$$
X^i=\begin{bmatrix}X_{11}^i&X_{12}^i\\(X_{12}^i)'&X_{22}^i\end{bmatrix},
$$
there are   $X_{22}^i=B'P^iB+D'P^iD+R$ and $X_{12}^i=A'P^iB+C'P^iD$,
which imply  that the primal update (\ref{primalupw}) is identical
to the PI  step (\ref{offlineF}). Hence, the dual and primal update
rules in Algorithm 1 can be interpreted as  a  policy evaluation and
a  policy improvement  in the off-line PI algorithm.  According to
Lemma \ref{converoff}, the  convergence property  can be directly
obtained.  $\blacksquare$

\begin{remark}
From the proof of Theorem \ref{thddddddd},   it is apparent that the
model-based primal-dual algorithm is equivalent to the  off-line  PI algorithm.
Thus,  the dual variable $X^i$ converges to the optimal
Q-function. Theorem~\ref{thddddddd} reveals the relation among the primal-dual algorithm, off-line PI and  Q-function.
\end{remark}
Notice that the matrix $\widetilde{S}$ can be estimated based on the
observations of the  state  and input variables. Next, we explore
the partially   model-free implementation of {Algorithm 2}, i.e.,
$A$ and $B$ are unknown, but $C$ and $D$ are known.  By limiting the
time frame of $\widetilde{S}$ defined in (\ref{sdefi2}), we
construct two new matrices
 \begin{eqnarray*}
\widetilde{S}(F^i) := \sum\limits_{l=1}^r\sum\limits_{k=0}^M
\mathcal{E}[v_k^{F^i, v^l}(v_k^{F^i, v^l})']
\end{eqnarray*}
and
\begin{eqnarray*}
W(F^i):= \sum\limits_{l=1}^r\sum\limits_{k=0}^{M}
\mathcal{E}[v_k^{F^i, v_l}(v_{k+1}^{F^i, v_l})']
=\widetilde{S}(F^i)A_F' .
\end{eqnarray*}
Then, due to the positive definiteness of $\widetilde{S}(F^i)$,
 the solvability  of  the dual update step (\ref{dualup}) is equivalently transformed into
 solving
\begin{eqnarray} \label{iiiddd}
&&\widetilde{S}(F^i)[  (A_{F^i})'X^iA_{F^i}+(C_{F^i})'X^iC_{F^i}+\Lambda-X^i]
\widetilde{S}(F^i)=0.
\end{eqnarray}
Therefore, (\ref{iiiddd}) holds iff
\begin{eqnarray}\label{freeal}
&&W(F^i)X^iW'(F^i)+\widetilde{S}(F^i)(C_{F^i})'X^iC_{F^i}
\widetilde{S}(F^i)+\widetilde{S}(F^i)(\Lambda-X^i)\widetilde{S}(F^i)=0.
\end{eqnarray}
From the analysis above, the partially model-free primal-dual algorithm is obtained.
\begin{algorithm}
  \caption{Partially Model-Free Primal-Dual Algorithm}
  \begin{algorithmic}[1]  
          \State { Initialization:} $F^0\in\mathscr{F}$,  the convergence tolerance
          $\varepsilon >0$, and the initial iteration  $i=0$;
           \State  { Repeat;}
            \State  Derive $\widetilde{S}(F^i)$ and $W(F^i)$ by calculating the mean-value $\mathcal{E}[v_k^{F^i, v^l}(v_k^{F^i, v^l})']$  based on $H$ sample
            paths {\bf  $\mathcal{V}_k^h$}:
           $$
            \mathcal{E}[v_k^{F^i, v^l}(v_k^{F^i, v^l})']\approx
            \frac{1}{H}\sum\limits_{h=1}^H[v_{k,h}^{F^i, v^l} (v_{k,h}^{F^i, v^l})'];
           $$
            \State Calculate $X^i$ by solving (\ref{freeal});
             \State  Update control gain as
               \begin{equation}\label{primalup}
             F^{i+1}=-(X_{22}^i)^{-1}(X_{12}^i)';
             \end{equation}
             \State \ \ \ \ \ \ \ $i\leftarrow i+1$;
             \State { Until} $\|F^i-F^{i+1}\|\leq \varepsilon$.
    \end{algorithmic}
\end{algorithm}

 \section{Simulation}
This section provides an example to evaluate the effectiveness of our
obtained results.  Consider  the sensorimotor control tasks studied
in \cite{shuzhan,pangboBC,shiyan}, where the human arm makes
horizontal point-to-point reach movements. The dynamics are
described by
\begin{eqnarray*}
dp=vdt, \ mdv=(a-bv+f)dt, \ \tau da=(u-a)dt+d\eta,
\end{eqnarray*}
where $p=[p_x\ p_y]'$, $v=[v_x\ v_y]'$, $a=[a_x\ a_y]'$, and $u=[u_x\ u_y]'$ represent the two-dimensional hand position,
velocity, actuator state, and control input, respectively. Notice that the term $f$ is used to model the external disturbance and we only consider
$f$  in the velocity-dependent force field
\begin{equation*}
f=\left[\begin{array}{cc}f_x\\ f_y\end{array}\right]
=\chi \left[\begin{array}{cc}13&-18\\ 18&13\end{array}\right]v
\end{equation*}
 with $\chi\in [2/3, 1]$, serving as an adjustable parameter based on the subject's strength.
 The term $d\eta$ is the control-dependent noise
  shown as
  $$
  d\eta=D_1ud\eta_1+D_2ud\eta_2
  $$
 with
 \begin{eqnarray*}
D_1=\left[\begin{array}{cc}d_1&0\\ d_2&0\end{array}\right], \
D_2=\left[\begin{array}{cc}0&-d_2\\0& d_1\end{array}\right].
\end{eqnarray*}
 $\eta_1$  and  $\eta_2$ are independent standard Wiener processes. The meaning and values of  parameters $m, b, \tau, c_1$  and $c_2$
 can be found in \cite{pangboBC}. In order to rewrite  this model in the state-space form,
 we take  the system state $x=[p\ v\ a]'$ and the discrete-time dynamic model is obtained by the  Euler discretization method with step size $\Delta t=0.1s$:
\begin{equation*}
x_{k+1}=Ax_k+Bu_k+D_1u_{k}w^1_k+D_2u_{k}w^2_k,
\end{equation*}
where
\begin{eqnarray*}
&&A=\left[\begin{array}{ccc}
0_2&I_2&0_2\\ 0_2&-\frac{b}{m}I_2&\frac{1}{m}I_2\\
0_2&0_2&-\frac{1}{\tau}I_2\end{array}\right], \
B=\left[\begin{array}{ccc}
0_2\\0_2\\\frac{1}{\tau}I_2\end{array}\right].
\end{eqnarray*}
The weighting matrices in cost function (\ref{cost1}) are given as
$Q=diag(\bar{Q}_1, \bar{Q}_2,\bar{Q}_3)$ and $R=0.01I_2$ with $\bar{Q}_3=diag(0.01, 0.01)$  and
\begin{eqnarray*}
\bar{Q}_1=\left[\begin{array}{cccccccc}
2000&-40 \\ -40&1000 \end{array}\right],
\
\bar{Q}_2=\left[\begin{array}{cccccccc}
20&-1\\
-1&20\end{array}\right].
\end{eqnarray*}
\begin{figure}[!htb]
  \centering
  \includegraphics[width=8.5cm]{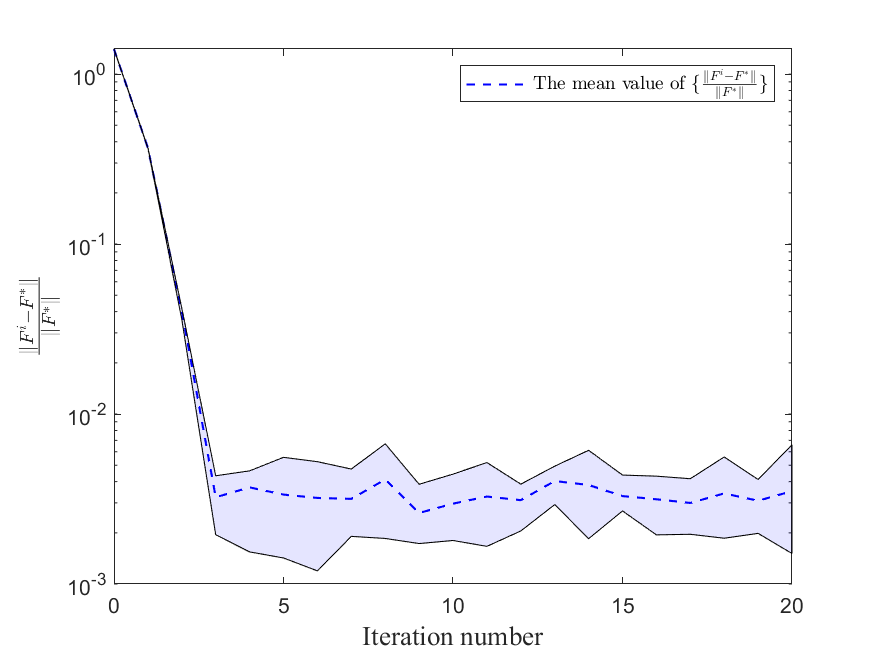}
  \caption{
The curves of
  relative error between the learned  control gain $F$
 and its true value $F^*$. }
  \label{hinfty}
\end{figure}
\begin{figure}[!htb]
  \centering
  \includegraphics[width=8.5cm]{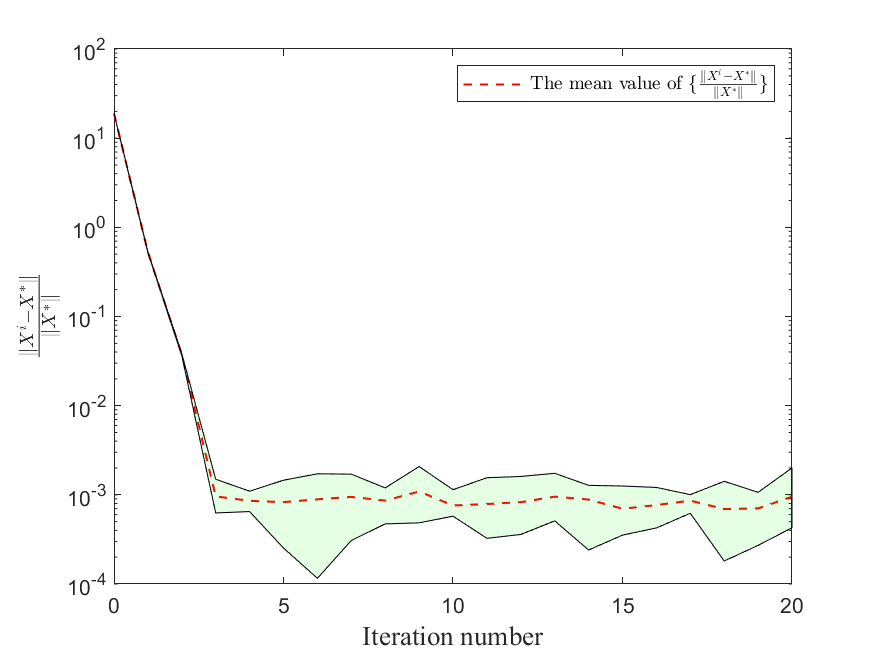}
  \caption{The curves of
  the relative error between the learned  Q-function $X$
 and its true value $X^*$.
 }
  \label{hinfty}
\end{figure}
\begin{figure}[!htb]
  \centering
  \includegraphics[width=8.5cm]{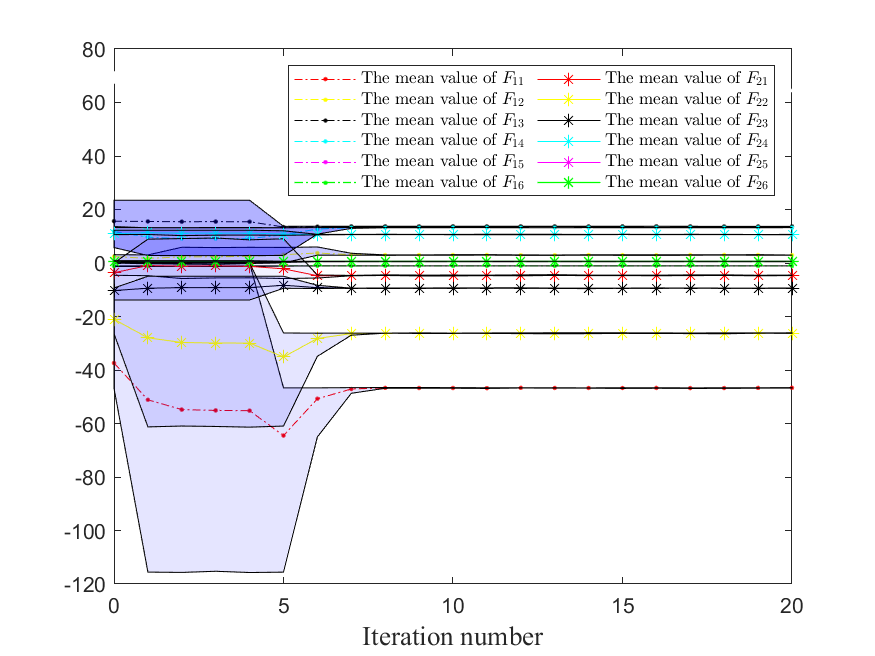}
  \caption{The curves of each element for learned control gain $F$ and its true value $F^*$. }
  \label{hinfty}
\end{figure}
Choosing the   initial  feedback  stabilizing  gain as
{\tiny\begin{equation*}
F^{(0)}=\left[\begin{array}{cccccccc}
 -0.0273   & -0.0258&    23.4596 &  5.7615&   0.2648&   -1.2886\\
    0.0238 &   0.0055  & -13.8178  &  12.2552  & 0.5310  &  0.8847
\end{array}
\right],
\end{equation*}}
the learned optimal control gain is achieved after $7$ iterations, which is
{\tiny\begin{equation*}
F^{(7)}=\left[\begin{array}{cccccccc}
-46.8727  &  3.0546  & 13.6769 &  13.2135 &   0.5656  & -1.0725\\
   -4.5767 & -26.8887 &  -9.3505&   10.6381 &   0.4850&    0.6393
\end{array}
\right].
\end{equation*}}
The Monte-Carlo experiment is adopted by implementing Algorithm 2 for $H=15$ times. We  terminate each experiment with  $i=20$.
Figure 1 illustrates the convergence of the relative error between the learned  control gain $F^i$ at each step
 and its true value $F^*$, i.e.,  $\frac{\|F^i-F^*\|}{\|F^*\|}$ with
{\tiny\begin{equation*}
F^{*}=\left[\begin{array}{cccccccc}
 -46.7316  &  2.9776 &  13.6867 &  13.2318  &  0.5664&  -1.0728\\
   -4.5899 & -26.2364&   -9.4265 &  10.6473  &  0.4846  &  0.6437
\end{array}
\right]
\end{equation*}}
as derived by the modified model-based PI algorithm
(\ref{offlineP})-(\ref{offlineF}). The  modifications
 involve replacing  $(C+DF^{(i)})$, $(D'P^{(i)}D)$, and $(D'P^{(i)}C)$
 by
$\sum_{m=1}^2(C_m+D_mF^{(i)})$,
$\sum_{m=1}^2(D'_mP^{(i)}D_m)$, and $\sum_{m=1}^2
(D'_mP^{(i)}C_m)$, respectively.
Similarly, Figure 2 presents the relative error convergence
 curve of the Q-function $X$.
 The dotted  lines  of  Figures 1 and 2  indicate the mean,  and the shaded areas of  Figures 1 and 2  cover    $15$ experimental
 trajectories  that show   the convergence of the learned values $F^i$ and $X^i$  using  Algorithm 2.
 Additionally, the relative error
  precisions  of $F^i$ and $X^i$  reach  $10^{-2}$ and $10^{-3}$ by the
  fourth iteration, respectively.
  To provide further insight into the convergence process, Figure 3 shows that each component of $F^i$ converges to
  its corresponding optimal component value.

\section{Conclusion}
In this paper, the primal-dual optimization method has been applied to solve  the
SLQR  of linear discrete-time stochastic systems including  model-based and model-free controller designs.      By
skillfully  constructing appropriate matrices $S$ and $X$,
the original dynamic optimization problem has been equivalently transformed
 into the solvability of  two matrix-valued  equations. Specially,  an augmented system that combines the state information and control
 input information   has been designed with any initial control input rather than
 depending on the initial state to  ensure the strict  positive definiteness of the modified
   matrix $\widetilde{S}$.
The strong duality  of   the primal and dual  problems
  was obtained  by choosing the  proper  dual variable $X$, from which the model-based
  primal-dual algorithm for equivalently solving the  SLQR problem   has been derived.
 Moreover, we have proved that the constructed dual variable $X$
  converges to the Q-function,  which  lays a new theoretical
  foundation for the Q-learning algorithm in RL.  A possible
extension is to use the primal-dual frame
to research  the  the fully  or    partially model-free robust $H_\infty$ control \cite{jiangauto} and mixed
$H_2/H_\infty$ \cite{zhang-xie-chen} control problems of  stochastic systems. Of course, there still remains  a few unsolved questions for SLQR issue including the fully model-free SLQR design ($A$, $B$, $C$,  and $D$ are all unknown) and the model-free design for  indefinite SLQR. All these problems merit further study.


\begin{thebibliography}{99}



\bibitem{anderson}
Anderson, B. D. O., \& Moore, J. B. (1989). \textit{Optimal
Control-Linear Quadratic Methods}. Prentice-Hall, New York.


\bibitem{rlcon2}
     Bertsekas, D. P., \& Tsitsiklis, J. N. (1996).
    \textit{Neuro-Dynamic Programming}. Belmont, MA, USA: Athena Scientific.



\bibitem{Busoniu}
Busoniu, L.,  de Bruin, T., Tolic, D.,  Kober, J., \&   Palunko, I. (2018).
Reinforcement learning for control: Performance, stability, and deep approximators. {\em Annual Reviews in Control},
46, 8-28.



\bibitem{boyd2004}  Boyd, B.,  \& Vandenberghe, V.  (2004). {\em Convex
Optimization}. Cambridge  University Press.
\bibitem{pangbo1}
Cui, L., Pang, B., Krstic, M., \& Jiang, Z.-P. (2025). Learning-based adaptive optimal control of linear time-delay systems: A value iteration approach. \textit{Automatica}, \textit{171}, 111944.





\bibitem{Bouhtouri_1999}
El  Bouhtouri, A.,  Hinrichsen, D.,     \&  Pritchard, A. J. (1999).
$H_\infty$-type control for discrete-time stochastic systems.  {\em
Int.  J.  Robust Nonlinear Control}, 9,  923-948.

\bibitem{che_98} Chen, S.,  Li, X.,  \&  Zhou, X. (1998). Stochastic linear
quadratic regulators with indefinite control weight costs.  {\em SIAM
J. Contr. Optim.,}  36, 1685-1702.



\bibitem{dom} Dombrovskii, V. V. \& Lyashenko,  E. A. (2003).  A linear quadratic control
for discrete systems with random parameters and multiplicative noise
and its application to investment portfolio optimization. {\em
Automat. Remote Control}, 64,  1558-1570.





\bibitem{farja}
    Farjadnasab, M., \& Babazadeh, M. (2022). Model-free LQR design by Q-function learning. \textit{Automatica}, \textit{137}, 110060.



\bibitem{Fazel}
Fazel, M., Ge, R., Kakade, S. M., \& Mesbahi, M. (2018). Global convergence of
policy gradient methods for the linear quadratic regulator. {\em In International
conference on machine learning},  1467-1476. Stockholm, Sweden.



    \bibitem{huangyulin}
    Huang, Y., Zhang, W., \& Zhang, H. (2008). Infinite horizon linear quadratic optimal control for discrete time stochastic systems.
    \textit{Asian Journal of Control}, \textit{10}(5), 608-615.


 \bibitem{jiangsci}
Jiang, X.,  Wang, Y.,  Zhao,  D., \& Shi, L. (2024).  Online Pareto optimal control of mean-field stochastic multi-player systems using policy iteration.   {\em Science China Information Sciences},  67(4),  140202:1-140202:17.


\bibitem{klein}
Kleinman, D. L. (1968), On an iterative technique
for Riccati equation computations, {\em IEEE Transactions on Automatic Control}, 13, 114-115.






    \bibitem{kalman}
    Kalman, R. E. (1960). Contributions to the theory of optimal control. \textit{Bol. Soc. Mat. Mex.}, \textit{5}(2), 102-119.


    \bibitem{Karl}
    Karl, K.  \&  Tu, S. (2019). Finite-time analysis of approximate policy iteration for
the linear quadratic regulator.   {\em In International conference on machine learning}, 8514-8524. Vancouver, Canada.


    \bibitem{jiangzp}
Kiumarsi, B., Lewis, F. L., \& Jiang, Z. P. (2017). $H_\infty $ control of linear discrete-time systems: Off-policy
    reinforcement learning. \textit{Automatica}, \textit{78}, 144-152.


\bibitem{lewis} Lewis, F. L. (1986). {\em Optimal Control.} John
 Wiley \& Sons.


    \bibitem{lee2019}
    Lee, D., \& Hu, J. (2019). Primal-dual Q-learning framework for LQR design. 
    \textit{IEEE Transactions on Automatic Control},
    \textit{64}(9), 3756-3763.

\bibitem{shuzhan}
    Lai, J., Xiong, J., \& Shu, Z. (2023). Model-free optimal control of discrete-time systems with additive and multiplicative noises. \textit{Automatica}, \textit{147}, 110685.



\bibitem{liman}
    Li, M., Qin, J., Zheng, W. X., Wang, Y., \& Kang, Y. (2022). Model-free design of stochastic LQR controller from a primal-dual
    optimization perspective. \textit{Automatica}, \textit{140}, 110253.





\bibitem{shiyan}
Liu, D.,  \& Todorov, E. (2007). Evidence for the flexible sensorimotor strategies predicted by optimal feedback control.
\textit{The Journal of Neuroscience}, \textit{27}(35), 93354-9368.





\bibitem{oura2024}
Oura, R., Ushio, T., \& Sakakibara, A. (2024). Bounded synthesis and reinforcement learning of supervisors for stochastic discrete event systems with LTL specifications. \textit{IEEE Transactions on Automatic Control}, \textit{69}(10), 6668--6683.


    \bibitem{pangbo}
    Pang, B., \& Jiang, Z. P. (2022). Reinforcement learning for adaptive optimal stationary control of linear stochastic systems.
    \textit{IEEE Transactions on Automatic Control}, \textit{68}(4), 2383-2390.



\bibitem{pangboBC}
Pang, B., Cui, L., \& Jiang, Z. P. (2022). Human motor learning is robust to control-dependent noise. \textit{Biological Cybernetics}, \textit{116}, 307-325.










\bibitem{sun}

Sun, J. \& J. Yong. (2023). Stochastic linear-quadratic optimal control problems-Recent
developments. {\em Annual Reviews in Control}, 56, 100899.

\bibitem{Stephen}
Stephen, T.  \&  Benjamin, R.  (2019). The gap between model-based
and model-free methods on the linear quadratic regulator: An asymptotic
 viewpoint. {\em In Proceedings of the Conference on Learning Theory},
3036-3083.



\bibitem{vrabie}
    Vrabie, D., Pastravanu, O., Abu-Khalaf, M., \& Lewis, F. L. (2009).
    Adaptive optimal control for continuous-time linear systems based on policy iteration. \textit{Automatica}, \textit{45}(2), 477-484.



\bibitem{wang2025}
Wang, Y., You, K., Huang, D., \& Shang, C. (2025). Data-driven output prediction and control of stochastic systems: An innovation-based approach. \textit{Automatica}, \textit{171}, 111897.




    \bibitem{watkins}
    Watkins, C. J., \& Dayan, P. (1992). Q-Learning. \textit{Machine Learning}, \textit{8}, 279-292.


\bibitem{rlcon1}
    Werbos, P. J. (1991). \textit{A Menu of Designs for Reinforcement Learning Over Time}. Cambridge, MA, USA: MIT Press.


\bibitem{wonham}
Wonham, W. M. (1968). On a matrix Riccati equation of stochastic
control.
    \textit{SIAM Journal on Control}, \textit{6}(4), 681-697.
\bibitem{ydd}
    Yao, D. D., Zhang, S., \& Zhou, X. Y. (2001). Stochastic linear-quadratic control via semidefinite programming.
    \textit{SIAM Journal on Control and Optimization}, \textit{40}(3), 801-823.


    \bibitem{ybbasi}
Yasin, A. Y., Nevena, L.,  \&  Csaba, S. (2019). Model-free
linear quadratic control via reduction to expert prediction. {\em In Proceedings
of the International Conference on Artificial Intelligence and Statistics},
3108-3117.

\bibitem{lixiuxian}
Yuan, K.,  Xu, W.,   \&  Ling, Q. (2020).  Can primal methods outperform primal-dual methods in decentralized dynamic optimization?.  {\em IEEE Transactions on Signal
Processing}, 68, 4466-4480.


\bibitem{zhanghan1}
Zhang, H., \& Ringh, A. (2023). Inverse linear-quadratic
discrete-time finite horizon optimal control for indistinguishable
homogeneous agents: A convex optimization approach.
\textit{Automatica},
    \textit{148}, 110758.



\bibitem{zwh2004}
    Zhang, W., \& Chen, B. S. (2004). On stabilizability and exact observability of stochastic systems with their
    applications. \textit{Automatica}, \textit{40}(1), 87-94.



\bibitem{zhang-xie-chen} Zhang, W.,  Xie, L., \& Chen, B. S.  (2017). {\em Stochastic  $H_2/H_\infty$  Control: A Nash Game
Approach}.   CRC Press.



\bibitem{zhang2008}
    Zhang, W., Zhang, H., \& Chen, B. S. (2008). Generalized Lyapunov equation approach to state-dependent stochastic
    stabilization/detectability criterion. \textit{IEEE Transactions on Automatic Control}, \textit{53}(7), 1630-1642.





\bibitem{zhangphd} Zhang, W. (1998). Study on the algebraic Riccati equation  arising from infinite horizon stochastic LQ
optimal control.  Zhejiang University, PhD dissertation.





\bibitem{jiangauto}
Zhang, W., Guo, J., \& Jiang, X. (2025).  Model-free $H_\infty$ control of It\^o stochastic system via off-policy reinforcement learning.  \textit{Automatica},  \textit{174},  112144.

\bibitem{jiangtac}
Zhang, W., Yu,  Z.,  \& Jiang, X. (2024).  Finite-time annular domain stability and asynchronous $H_\infty$  control for stochastic switching Markov jump systems.  \textit{IEEE Transactions on Automatic Control}, \textit{69}{9}, 6277-6284.
\bibitem{youbased}
Zhao, B., \& You, K. (2023). Survey of recent progress in data-driven policy optimization for controller design (in Chinese). \textit{Scientia Sinica Informationis}, \textit{53}(6), 1027--1049.
\end{thebibliography}
\end{document}